
\def\Z{\gz}
\def\im{\mathop{\rm Im}\nolimits} 
\def\semidirect{\mathop{:}\nolimits} 
\def\isomorphism{\cong}  
\def\narukismodel{\check{C}}
\def\smallmatrix#1#2#3#4{%
	\left({#1\atop #3}\;{#2\atop #4}\right)}
\def\Eisenstein{\calE}

\font\pro =cmss10
\catcode`\^^Z=9
\catcode`\^^M=10
\output={\if N\header\headline={\hfill}\fi
\plainoutput\global\let\header=Y}
\magnification\magstep1
\tolerance = 500
\hsize=14.4true cm
\vsize=22.5true cm
\parindent=6true mm\overfullrule=2pt
\newcount\kapnum \kapnum=0
\newcount\parnum \parnum=0
\newcount\procnum \procnum=0
\newcount\nicknum \nicknum=1
\font\ninett=cmtt9

\font\ninebf=cmbx9

\font\sixbf=cmbx6
\font\ninesl=cmsl9

\font\nineit=cmti9

\font\ninerm=cmr9

\font\sixrm=cmr6
\font\ninei=cmmi9
\font\eighti=cmmi8
\font\sixi=cmmi6
\skewchar\ninei='177 \skewchar\eighti='177 \skewchar\sixi='177
\font\ninesy=cmsy9
\font\eightsy=cmsy8
\font\sixsy=cmsy6
\skewchar\ninesy='60 \skewchar\eightsy='60 \skewchar\sixsy='60
\font\titelfont=cmr10 scaled 1440
\font\paragratit=cmbx10 scaled 1200

\font\name=cmcsc10
\font\emph=cmbxti10

\font\tenmsbm=msbm10
\font\sevenmsbm=msbm7
%

%
\font\got=eufm10

\font\teneufm=eufm10
\font\seveneufm=eufm7
\font\fiveeufm=eufm5
\newfam\eufmfam
\textfont\eufmfam=\teneufm
\scriptfont\eufmfam=\seveneufm
\scriptscriptfont\eufmfam=\fiveeufm

\font\tenmsam=msam10
\font\sevenmsam=msam7
\font\fivemsam=msam5
\newfam\msamfam
\textfont\msamfam=\tenmsam
\scriptfont\msamfam=\sevenmsam
\scriptscriptfont\msamfam=\fivemsam
\font\tenmsbm=msbm10
\font\sevenmsbm=msbm7
\font\fivemsbm=msbm5
\newfam\msbmfam
\textfont\msbmfam=\tenmsbm
\scriptfont\msbmfam=\sevenmsbm
\scriptscriptfont\msbmfam=\fivemsbm
\def\Bbb#1{{\fam\msbmfam\relax#1}}
\def\cz{{\kern0.4pt\Bbb C\kern0.7pt}
}
\def\ez{{\kern0.4pt\Bbb E\kern0.7pt}
}
\def\fz{{\kern0.4pt\Bbb F\kern0.7pt}}
\def\gz{{\kern0.4pt\Bbb Z\kern0.7pt}}
\def\hz{{\kern0.4pt\Bbb H\kern0.7pt}
}
\def\kz{{\kern0.4pt\Bbb K\kern0.7pt}
}
\def\nz{{\kern0.4pt\Bbb N\kern0.7pt}
}
\def\rz{{\kern0.4pt\Bbb R\kern0.7pt}
}
\def\sz{{\kern0.4pt\Bbb S\kern0.7pt}
}
\def\pz{{\kern0.4pt\Bbb P\kern0.7pt}
}
\def\qz{{\kern0.4pt\Bbb Q\kern0.7pt}
}
\newskip\ttglue
\def\ninepoint{\def\rm{\fam0\ninerm}%
  \textfont0=\ninerm \scriptfont0=\sixrm \scriptscriptfont0=\fiverm
  \textfont1=\ninei \scriptfont1=\sixi \scriptscriptfont1=\fivei
  \textfont2=\ninesy \scriptfont2=\sixsy \scriptscriptfont2=\fivesy
  \textfont3=\tenex \scriptfont3=\tenex \scriptscriptfont3=\tenex
  \def\it{\fam\itfam\nineit}%
  \textfont\itfam=\nineit
  \def\sl{\fam\slfam\ninesl}%
  \textfont\slfam=\ninesl
  \def\bf{\fam\bffam\ninebf}%
  \textfont\bffam=\ninebf \scriptfont\bffam=\sixbf
   \scriptscriptfont\bffam=\fivebf
  \def\tt{\fam\ttfam\ninett}%
  \textfont\ttfam=\ninett
  \tt \ttglue=.5em plus.25em minus.15em
  \normalbaselineskip=11pt
  \font\name=cmcsc9
  \let\sc=\sevenrm
  \let\big=\ninebig
  \setbox\strutbox=\hbox{\vrule height8pt depth3pt width0pt}%
  \normalbaselines\rm
  \def\sl{\it}}

\headline={\ifodd\pageno\rightheadline\else\leftheadline\fi}
\def\rightheadline{\ninepoint Paragraphen"uberschrift\hfill\folio}
\def\leftheadline{\ninepoint\folio\hfill Chapter"uberschrift}
\let\header=Y
\def\titel#1{\need 9cm \vskip 2truecm
\parnum=0\global\advance \kapnum by 1
{\baselineskip=16pt\lineskip=16pt\rightskip0pt
plus4em\spaceskip.3333em\xspaceskip.5em\pretolerance=10000\noindent
\titelfont Chapter \uppercase\expandafter{\romannumeral\kapnum}.
#1\vskip2true cm}\def\leftheadline{\ninepoint
\folio\hfill Chapter \uppercase\expandafter{\romannumeral\kapnum}.
#1}\let\header=N
}
\def\Titel#1{\need 9cm \vskip 2truecm
\global\advance \kapnum by 1
{\baselineskip=16pt\lineskip=16pt\rightskip0pt
plus4em\spaceskip.3333em\xspaceskip.5em\pretolerance=10000\noindent
\titelfont\uppercase\expandafter{\romannumeral\kapnum}.
#1\vskip2true cm}\def\leftheadline{\ninepoint
\folio\hfill\uppercase\expandafter{\romannumeral\kapnum}.
#1}\let\header=N
}
\def\need#1cm {\par\dimen0=\pagetotal\ifdim\dimen0<\vsize
\global\advance\dimen0by#1 true cm
\ifdim\dimen0>\vsize\vfil\eject\noindent\fi\fi}
\def\neupara#1{\par\penalty-2000
\procnum=0\global\advance\parnum by 1
\vskip1cm\noindent{\paragratit \the\parnum. #1}%
\def\rightheadline{\ninepoint\S\the\parnum.\ #1\hfill \folio}%
\vskip 8mm\noindent}
\def\Proclaim #1 #2\finishproclaim {\bigbreak\noindent
{\bf#1\unskip{}. }{\it#2}\medbreak\noindent}
%
\gdef\proclaim #1 #2 #3\finishproclaim {\bigbreak\noindent%
\global\advance\procnum by 1 
{%
{\relax\ifodd \nicknum 
\hbox to 0pt{\vrule depth 0pt height0pt width\hsize
   \quad \ninett#3\hss}\else {}\fi}%
\bf\the\parnum.\the\procnum\ #1\unskip{}. }
{\it#2}
\medbreak\noindent}
\newcount\stunde \newcount\minute \newcount\hilfsvar
\def\uhrzeit{
	\stunde=\the\time \divide \stunde by 60
	\minute=\the\time 
	\hilfsvar=\stunde \multiply \hilfsvar by 60
	\advance \minute by -\hilfsvar
	\ifnum\the\stunde<10
	\ifnum\the\minute<10
	0\the\stunde:0\the\minute~Uhr
	\else
	0\the\stunde:\the\minute~Uhr
	\fi
	\else
	\ifnum\the\minute<10
	\the\stunde:0\the\minute~Uhr
	\else
	\the\stunde:\the\minute~Uhr
	\fi
	\fi
	}
 \def\calB{{\cal B}}
     
\def\calC{{\cal C}} 
\def\calE{{\cal E}} 
 \def\calH{{\cal H}}

\def\calM{{\cal M}}

\def\calU{{\cal U}} \def\calV{{\cal V}}

\def\gota{\hbox{\got a}}

\def\goto{\hbox{\got o}}

\def\schreib#1{\hbox{#1}}

\def\Aut{\mathop{\rm Aut}\nolimits}

\def\im{\mathop{\rm Im}\nolimits} \def\Im{\im}

\def\mod{\mathop{\rm mod}\nolimits}
\def\O{{\rm O}} 
\def\U{{\rm U}}

\def\SL{\mathop{\rm SL}\nolimits}

\def\boxit#1{
  \vbox{\hrule\hbox{\vrule\kern6pt
  \vbox{\kern8pt#1\kern8pt}\kern6pt\vrule}\hrule}}
\def\Boxit#1{
  \vbox{\hrule\hbox{\vrule\kern2pt
  \vbox{\kern2pt#1\kern2pt}\kern2pt\vrule}\hrule}}

\def\zwischen#1{\bigbreak\noindent{\bf#1\medbreak\noindent}}

\def\smallni{\smallskip\noindent }
\def\medni{\medskip\noindent }

\def\lo{\longrightarrow}

\def\loma{\longmapsto}
\def\betr#1{\vert#1\vert}
\def\spitz#1{\langle#1\rangle}
\def\pii{\pi {\rm i}}

\def\mag{\hbox{\rm i}}
\def\square{\hbox{\hbox to 0pt{$\sqcup$\hss}\hbox{$\sqcap$}}}
\def\qed{\ifmmode\square\else{\unskip\nobreak\hfil
\penalty50\hskip3em\null\nobreak\hfil\square
\parfillskip=0pt\finalhyphendemerits=0\endgraf}\fi}
\def\pn{\the\parnum.\the\procnum}
\def\downmapsto{{\buildrel
        {\vbox{\hbox{\hskip.2pt$\scriptstyle-$}}}
        \over{\raise7pt\vbox{\vskip-4pt\hbox{$\textstyle\downarrow$}}}}}
\newcount\stunde \newcount\minute \newcount\hilfsvar
\def\uhrzeit{
	\stunde=\the\time \divide \stunde by 60
	\minute=\the\time 
	\hilfsvar=\stunde \multiply \hilfsvar by 60
	\advance \minute by -\hilfsvar
	\ifnum\the\stunde<10
	\ifnum\the\minute<10
	0\the\stunde:0\the\minute~Uhr
	\else
	0\the\stunde:\the\minute~Uhr
	\fi
	\else
	\ifnum\the\minute<10
	\the\stunde:0\the\minute~Uhr
	\else
	\the\stunde:\the\minute~Uhr
	\fi
	\fi
	}
\def\heute{
	\number\day.\ts
	\ifcase\month
	\or1\or2\or3\or4\or5\or6\or7\or8\or9\or10\or11\or12\fi
	.\ts\number\year
	}

\def\Orb{2.1}

\def\Fraemp{3.4}
\def\foo{3.5}
\def\Badd{4.1}
\def\constr{4.2}
\def\Sing{4.3}
\def\addEm{4.4}
\def\EinF{4.5}
\def\farD{4.6}
\def\farDI{4.7}
\def\dimeis{5.1}
\def\bopro{5.2}
\def\noobstr{5.3}

\def\daSi{5.6}
\def\BiR{6.1}

\def\cubrel {6.3}

\def\crdivisor{7.1}
\def\ratios{7.2}
\def\embedding{7.3}
%
\nicknum=0      
\def\RAND#1{\hbox to 0mm{\hss\vtop to 0pt{%
  \raggedright\ninepoint\parindent=0pt%
  \baselineskip=1pt\hsize=2cm #1\vss}}}

\def\frac#1#2{#1/#2}

\noindent{\titelfont Cubic surfaces and Borcherds products}%
\def\leftheadline{\ninepoint\folio\hfill 
    Cubic surfaces and Borcherds products}%
\def\leftheadline{\ninepoint \folio\hfill Cubic surfaces and 
     Borcherds products}%
\def\rightheadline{\ninepoint Introduction\hfill \folio}%
\headline={\ifodd\pageno\rightheadline\else\leftheadline\fi}
\def\tensor{\otimes}

\def\imag{\hbox{\rm i}}

\def\bz{\calB}
\vskip1cm\noindent
\centerline{\paragratit \rm February 7, 2000}
\centerline{2000 MSC: 11F55, 14J10}
\vskip5mm\noindent
\centerline{
\vbox{\openup-3pt\noindent\hsize=6cm{\it
		Daniel Allcock
\hfil\break Department of Mathematics
\hfil\break Harvard University
\hfil\break Cambridge, MA 02138
}\hfil\break
{\sevenrm
ALLCOCK@\vskip-1mm\noindent
\null\quad MATH.HARVARD.EDU}}
\hfill
\vbox{\openup-3pt\noindent\hsize=6cm{\it Eberhard Freitag
\hfil\break Mathematisches Institut
\hfil\break Im Neuenheimer Feld 288
\hfil\break D69120 Heidelberg
\hfil\hfil\break}{\sevenrm
FREITAG@\vskip-1mm\noindent
\null\quad MATHI.UNI-HEIDELBERG.DE}}}
\let\header=N%
\medskip\neupara{Introduction}%
The moduli space $\calM$ of 
marked cubic surfaces can be identified with
the Baily-Borel compactification of $\bz_4/\Gamma$,
where $\bz_4$ denotes the complex $4$-ball and $\Gamma$ is a certain
arithmetic reflection group. (See [ACT2] and also [ACT1].)
In this paper we use the methods of R.~Borcherds to construct
automorphic forms on $\bz_4$. We will
obtain an embedding of $\calM$
into the $9$-dimensional projective space $P^9(\cz)$, whose
image is the intersection of 
270 explicitly known cubic
$8$-folds.
This map is compatible with the actions of the Weyl
group $W(E_6)$ on $\calM$ and $P^9$. The former action arises
because
$W(E_6)$ permutes the markings of cubic surfaces, and the 
latter action arises from the unique irreducible $10$-dimensional
representation of $W(E_6)$. Furthermore, the cubic $8$-folds 
are all equivalent under  $W(E_6)$.
\smallskip
The $10$-dimensional linear system associated to this map into
$P^9(\cz)$ contains $270$ automorphic forms with known zeros, which play a
central role in our investigation. 
In particular, there is a direct connection between them and the
classical invariants of cubic surfaces introduced by Cayley. He
considered the 27 lines on a smooth cubic surface and a certain
configuration of 45 planes that they determine. By considering
4-tuples of these planes that meet along one of the 27
lines, Cayley constructed 270 cross-ratios, and showed that these
allow one to recover the original surface. We show that Cayley's
cross-ratios coincide not with our Borcherds products but rather
with the quotients of certain pairs of them. This relies
on work of Naruki [Na] and is the main part of our proof that
our map of $\cal M$ into $P^9$ is an embedding.
\smallskip
We are grateful to R.~Borcherds, B.~van~Geemen, and R.~Vakil
for helpful discussions.
\neupara{The complex reflection group}%
Let
$$\calE=\gz[\omega],\quad
  \omega=[\root 3 \of 1]=-{1\over2}+{\sqrt{-3}\over2},$$
be the ring of Eisenstein integers.
We consider the lattice
$${\Lambda}=\calE^{1,4},$$
which is the $\calE$-module $\calE^5$
equipped with the hermitian form of signature $(1,4)$ given by
$$
\spitz{a,b}=\bar a_0b_0-\bar a_1b_1-\cdots-\bar a_4b_4.
\eqno (2.1)
$$
Let $\Aut({\Lambda})$ be the 
unitary group of this lattice, i.e.\ the
group of $\calE$-module automorphisms which preserve the
hermitian form. 
Complex conjugation acts as the identity on the residue field
$$\fz_3=\calE/\sqrt{-3}\,\calE,$$
which has order 3,
so the hermitian form induces a $\fz_3$-valued
quadratic form on the 5-dimensional $\fz_3$-vector space
${V}={\Lambda}/\sqrt{-3}\,{\Lambda}$. We denote the orthogonal
group of ${V}$ by $\O(5,3)$ and define $\Gamma$ to 
be the kernel of the action of $\Aut({\Lambda})$
on ${V}$.
We have the exact sequence
$$1\lo\Gamma\lo\Aut({\Lambda})\lo\O(5,3)\lo1.$$
For future reference we mention that 
${V}$ contains 242 nonzero elements, of
which 80 have norm $0$, 90 have norm $1$ and 72 have norm $-1$.
Nonzero vectors in ${V}$ are equivalent 
under $\O(5,3)$ if and only if they have
the same norm.
The subgroup of 
$\O(5,3)$  generated by 
the reflections in the norm $-1$ vectors is isomorphic to the
Weyl group $W(E_6)$, and $\O(5,3)\isomorphism
W(E_6)\times\{\pm1\}$. Furthermore,
$W(E_6)$ contains a simple subgroup of index~2 and 
order $25\,920$. 
\smallskip
A lattice vector $a\in {\Lambda}$ is called
primitive if it cannot be divided in ${\Lambda}$ by a non-unit of $\calE$.
Also, $a$ is called
$$
\vbox{\halign{\qquad\it#\hfil&\quad#\hfil\cr
isotropic  &if $\;\spitz{a,a}=0$,\cr
a short root &if
  $\;\spitz{a,a}=-1$, or\cr
a long root &if
  $\;\spitz{a,a}=-2$.\cr
}}
$$
The roots are important because $\Aut({\Lambda})$ contains reflections
in them.
If $a$ is a short root and $\zeta$ is a unit of $\calE$ (a sixth
root of unity) then the map
$$v\loma v-(1-\zeta){\spitz{a,v}\over\spitz{a,a}}a$$
is an automorphism of ${\Lambda}$. 
(In the special case $\zeta=\pm1$
this is also true if $a$ is a long root.)
This automorphism fixes the orthogonal complement of $a$ and maps
$a$ to $\zeta a$.
We call this automorphism a reflection if $\zeta\ne 1$.
The order of a reflection is two, three or six, and
we sometimes call reflections of these orders biflections,
triflections and hexflections. 
The
third roots of unity are congruent to 1 mod $\sqrt{-3}$, and
therefore the triflections belong to the congruence group
$\Gamma$. We remark that these triflections actually generate
$\Gamma$ [ACT2], although we will not need this fact.
\smallskip
We need some information about the orbit structure of ${\Lambda}$
with respect to $\Gamma$. 
If $a,b\in {\Lambda}$ are in the same $\Gamma$-orbit, then their
images in ${V}$ coincide. In some special cases the converse
is true:
\proclaim
{Proposition}
{Let $a$ and $b$ be 
two primitive
isotropic vectors, or two short roots, or
two long roots.
Then 
$a$ and $b$ are equivalent under $\Gamma$ if and only if their
images in ${V}$ coincide. The number of\/ $\Gamma$-orbits
of lines $\cz a$, where $a$ is a primitive isotropic vector,
a short root or a long root, is
40, 36 or 45, respectively. } 
Orb%
\finishproclaim
{\it Proof.}
The ``only if'' part is trivial.
To prove the converse,
we
use the fact that $\Aut({\Lambda})$ acts transitively on primitive
isotropic vectors, on short roots, and on long roots (see
Theorems~7.22 and~9.15 of [ACT2]).
It is a general fact that if a group $G$ acts transitively on a
set $X$, $N$ is a normal subgroup, and $x\in X$ has stabilizer
$G_x$ in $G$, then the orbits of $N$ on $X$ are in 1-1
correspondence with the cosets in $G/N$ of the image of
$G_x$. We apply this with $G=\Aut({\Lambda})$, $N=\Gamma$, and $x$
a primitive isotropic vector, short root or long
root of ${\Lambda}$. Then the number of $\Gamma$-orbits into which the
$\Aut({\Lambda})$-orbit of $x$ splits is equal to the index in $\O(5,3)$
of the reduction modulo $\sqrt{-3}$ of $\Aut({\Lambda})_x$. We will now
compute these reductions.
\smallskip
We first take $x$ to be a primitive null vector. 
According to paragraph~7.8 of [ACT2], its stabilizer in $\Aut({\Lambda})$
contains as a normal subgroup a Heisenberg group with center
$\im(\calE)$ and central quotient $\calE^3$, and the
stabilizer modulo this Heisenberg group is the isometry group
$(\gz/6)^3\semidirect S_3$ of the lattice $\calE^3$. By
considering the matrices for these transformations, it is easy
to see that the center of the Heisenberg group acts trivially on
${V}$, that $\calE^3$ acts as $\calE^3/(\sqrt{-3}\,\calE^3)\isomorphism
(\gz/3)^3$, that $(\gz/6)^3$ acts as $(\gz/2)^3$, and that
$S_3$ acts faithfully. The image of the stabilizer in
$\O(5,3)$ is a group $3^3\semidirect2^3\semidirect S_3$, which
has index 80 in $\O(5,3)$. Next we take $x$ to be a short root of
${\Lambda}$, say $(0,0,0,0,1)$, and $\bar x$ to be its image in
${V}$. Then the stabilizer of $\bar{x}$ is the 
orthogonal group of $\bar{x}^\bot$, which is
generated by the reflections in the nonisotropic elements of
$\bar{x}^\bot$. One can enumerate these vectors and check that
each is the image of a  root of $x^\bot$. The
biflections in these roots reduce to reflections of
${V}$, proving that the stabilizer of $x$ in $\Aut({\Lambda})$ surjects to the
stabilizer of $\bar{x}$, which has index 72 in $\O(5,3)$.
Exactly the same argument applies if
$x$ is a long root, say $(0,0,0,1,1)$, yielding an index of 90.
\smallskip
We have shown that there are 80 (resp. 72, 90) orbits of
primitive isotropic vectors (resp. short roots, long roots) in
${\Lambda}$, which is the same as the number of nonzero elements of
${V}$ of norm 0 (resp. $-1$, $1$). Since the map from
$\Gamma$-orbits of such lattice vectors to the corresponding set
of vectors in $V$ is onto, it is bijective.
This proves the first claim of the theorem, and the second
follows immediately.
\qed
\neupara{The ball quotient}%
The group $\Gamma$ acts on a complex 4-ball in the projective
space of $\cz^{1,4}=\Lambda\tensor_\calE\cz$. We will describe this in
some generality, for convenience in later sections.
Let $\goto$ be an order in an imaginary quadratic number
field. An $\goto$-lattice $L$ is a finitely generated
projective $\goto$-module
equipped with a Hermitian pairing $\spitz{\,,}$ on $L$ that takes
value in the field of fractions of $\goto$.
We take such pairings to be antilinear in the first and linear
in the second variable. We say that $L$ is Lorentzian when its
signature is $(1,n)$ with $n\geq1$. A point of the projective
space $P(L\tensor_{\goto}\cz)$ is called positive if it is
represented by a vector of positive norm. When $L$ is
Lorentzian, the positive points form an open $n$-ball $\calB(L)$
in projective space, which is
also called the complex hyperbolic space of $L$. 
$\Aut(L)$ acts properly discontinuously 
on $\calB(L)$, and there is a natural
compactification of the quotient, due to Baily and Borel [BB]. A
cusp is an element of $P(L\tensor_{\goto}\cz)$ that can be
represented by an isotropic lattice vector. The cusps are the
rational boundary points of $\calB(L)$, and there are only
finitely many orbits under $\Aut(L)$. We denote by
$\calB^*(L)$ the union of $\calB(L)$ with the set of all
cusps. The group $\Aut(L)$ acts on this extension, 
and the quotient of $\calB^*(L)$ by any finite-index
subgroup of $\Aut(L)$ carries the structure of a projective
algebraic variety.
\smallskip
In our setting we have $\goto=\calE$ and $L=\Lambda$. The
hermitian form on $\cz^{1,4}=\Lambda\tensor_\calE\cz$ is given by
Eq.~(2.1), 
and the identification of $\calB(\Lambda)$ with the
complex 4-ball
is easy. Namely, any element of $\calB(\Lambda)$ has a unique
representative $z\in\cz^{1,4}$ whose $z_0$-component is
$1$. Considering the remaining coordinates identifies
$\calB(\Lambda)$ with the set of all $(z_1,\ldots,z_4)\in\cz^4$
satisfying 
$$
\betr{z_1}^2+\cdots+\betr{z_4}^2<1.
\eqno (3.1)
$$
We will write $\calB_4$ for $\calB(\Lambda)$. We are interested
in the quotient $X$ of $\calB^*_4$ by $\Gamma\subset\Aut(\Lambda)$.
By Prop.~\Orb, there are forty $\Gamma$-orbits of
cusps in $\calB_4^*$, so the boundary of the 
Baily-Borel compactification of
$\calB_4/\Gamma$ consists of 40 points.
\smallskip
Let $a\in {\Lambda}$ be a vector of negative
norm.  The orthogonal complement $a^\perp$
of $a$ in $P(\cz^{1,4})$
meets $\calB_4$ nontrivially 
because $a$ has negative norm. We can consider
its intersection with $\calB_4^*$ and the image of this in $X$.
It is known that this is an algebraic subvariety of codimension
one, and
we are interested in this construction for $a$ a root of $\Lambda$.
In this case we call 
$a^\perp\cap\calB_4^*$
a mirror of $\calB_4^*$. 
The terminology derives from the fact that the mirror is the
fixed-point set of the reflection(s) in $a$, and we call
the mirror
short or long according to whether $a$ is short or long. 
The image in $X$ of a short (long) mirror
is called a short (long) mirror of $X$.
For convenience we sometimes call a vector in ${V}$ short
(resp. long)
if it has norm $-1$ (resp. 1). The short (long) vectors in ${V}$
are exactly the images of the short (long) roots of $\Lambda$.
The short (long) mirrors in $X$ correspond to the
36 (45) pairs $\{\pm a\}$ of short (long) vectors of ${V}$. 
We will need some results about the intersection behavior
of mirrors.
Orthogonality of mirrors in $\calB_4^*$ is defined in the obvious
way, and
we call two mirrors in $X$ orthogonal if the corresponding
elements of ${V}$ are orthogonal. 
If two mirrors in $\calB_4^*$
are orthogonal then so are their images in $X$.
\proclaim
{Lemma}
{Two short mirrors in $\calB_4$ are either orthogonal or
disjoint. 
}
emptzw%
\finishproclaim
{\it Proof}. 
We take $x$ and $y$ to be short roots whose mirrors are the
given mirrors. If the mirrors meet in $\calB^4$ then
a point of the intersection represents
a positive-definite one-dimensional subspace of $\cz^{1,4}$. Its
orthogonal complement is negative definite and contains $x$
and $y$. Hence the Gram matrix of $x,y$ must be positive, so
$\spitz{x,x}\spitz{y,y}-\betr{\spitz{x,y}}^2>0$, so
$\betr{\spitz{x,y}}^2<1$. Since $\spitz{x,y}\in\calE$ we must
have $\spitz{x,y}=0$.\qed
\smallskip
We now introduce the notion of a cross. This is fundamental 
for the paper because
the automorphic forms we will construct vanish exactly along the
points of a {cross} in $\calB_4^*$. The word ``cross'' is meant to
suggest several  mutually orthogonal objects.
\proclaim
{Definition}
{A {\emph {cross}} in ${V}$ is a set of 5 pairwise orthogonal
pairs $\pm a_i$, one pair consisting of long vectors and the
others consisting of short vectors.
The associated {\emph{cross}} in $X$ is the union of the
mirrors of the $\pm a_i$; it follows that a {cross} in
$X$ is a set of 5 pairwise orthogonal 
mirrors, one long and 4 short. The associated {\emph{cross}} in
$\calB_4^*$ is the preimage of the {cross} in $X$. A point
of $\calB_4^*$ lies in this {cross} just if it is orthogonal to a
root whose projection to ${V}$ is one of the $\pm a_i$.}
Frame%
\finishproclaim
Since
the three types of {cross} are in natural bijection, we will pass
between them without comment.
\proclaim{Lemma}
{There are 135 {crosses}, three containing each of the 45
long mirrors of $X$, and all 135 crosses are all equivalent
under $\O(5,3)$.
More precisely, if $\ell$ is a long mirror in $X$ then the 12 short mirrors
orthogonal to $\ell$ decompose
in a unique way 
into three sets of 4 mirrors which are pairwise orthogonal, and
the stabilizer of $\ell$
in $\O(5,3)$
permutes these sets transitively. 
}
deco%
\finishproclaim
{\it Proof.} 
The transitivity of $\O(5,3)$ on crosses in $V$ is
obvious, and the rest is just a calculation. Namely, the
orthogonal complement of a long vector $a$ contains 12 pairs
$\{\pm v\}$ of short vectors, orthogonality is (surprisingly) a
transitive relation, generating an equivalence relation with three classes of
size 4. By symmetry it suffices to check this for a single long
vector $a$, say $(1,0,0,0,0)$. Then the three classes are
$$
\vbox{%
\halign{$#$\ \hfil&$#$\ \hfil&$#$\ \hfil&$#$\hfil\cr
\{
 \pm(0,1,0,0,0),
&\pm(0,0,1,0,0),
&\pm(0,0,0,1,0),
&\pm(0,0,0,0,1)
\},\cr\{
 \pm(0,-1, 1, 1, 1),
&\pm(0, 1,-1, 1, 1),
&\pm(0, 1, 1,-1, 1),
&\pm(0, 1, 1, 1,-1)
\},\cr\{
 \pm(0,1, 1, 1, 1),
&\pm(0,1, 1,-1,-1),
&\pm(0,1,-1, 1,-1),
&\pm(0,1,-1,-1, 1)
\}.
\cr}}
$$
\qed
The purpose of the following theorem is to allow us to prove in
section~4 
that the automorphic forms we construct there have no
common zeros. 
\proclaim
{Theorem}
{No point of $X$ lies on all 135 {crosses}. Furthermore, if $p$
is the point of $\calB_4$ represented by
$(1,0,0,0,0)\in {\Lambda}$, then the image of
$p$ in $X$ is the only point of $X$ that lies on all the crosses
containing it.
Finally, for each boundary point $b$ of $X$, $b$ is the only
point of $X$ that lies on all the crosses containing it.}
Fraemp%
\finishproclaim
\smallskip
In order to prove the theorem we will need to understand the
$\Gamma$-orbits of points of $\calB_4$ that, like $p$, lie on
four short mirrors. If $q$ is such a point, then $q^\bot$ is a
copy of the unimodular lattice $\calE^{0,4}$, and it follows that
$q$ is represented by a lattice vector of norm $1$, and indeed
by six such vectors. The
images in ${V}$ of these vectors and of the short roots of
$q^\bot$ form a {cross}, which we call the {cross} associated to
$q$. 
\proclaim
{Lemma}
{The map just defined, which associates a {cross} to each point of
$\calB_4$ that lies on $4$ short mirrors, defines a bijection
between the set of $\Gamma$-orbits of such points and the set of
{crosses}. If each of $p,q\in \calB_4$ lies on four short mirrors,
and the images in ${V}$ of the short roots of $p^\bot$
coincide with the images of the short roots of $q^\bot$,
then $p$ and $q$ are $\Gamma$-equivalent.}
foo%
\finishproclaim
{\it Proof.}
For the first claim one uses the argument of \Orb. The essential
facts are that $\Aut({\Lambda})$ acts transitively on such
points of $\calB_4$ and that the stabilizer in $\Aut({\Lambda})$
of such a point of 
$\calB_4$ is $(\gz/6)\times(\gz/6)^4\semidirect S_4$, which
reduces modulo $\sqrt{-3}$ to $(\gz/2)^5\semidirect S_4$, of index
135 in $\O(5,3)$. 
(The transitivity statement
follows from the fact that such points in $\calB_4$
correspond bijectively to
the decompositions of $\Lambda$ as a direct sum
$\calE^{1,0}\oplus\calE^{0,4}$.)
The second claim is 
a consequence of the first: the short vectors of a {cross}
determine the {cross} uniquely, so the {crosses} associated to $p$
and $q$ coincide. \qed
\smallskip\noindent
{\it Proof of Theorem~\Fraemp.}
Most of the proof consists of computer calculations concerning
combinatorics in $V$; we will describe the ideas in sufficient detail for
them to be reproduced easily. 
\smallskip
One can enumerate the roots
orthogonal to $p$, and their images in ${V}$. A {cross} contains
$p$ just if it contains one of these images. One can compute the
set $\calC$ of {crosses} satisfying this condition, and one finds
$|\calC|=69<135$. In particular, $p$ does not lie on all~135
crosses. Now we will show that $p$ is the only point of
$\calB_4$ (up to
$\Gamma$-equivalence) that lies on all the crosses containing $p$.
Suppose $q\in\calB_4$ lies on every cross of $\calC$; we
will show that $q$ is $\Gamma$-equivalent to $p$.
First we will show that $q$ lies on 4 short mirrors. 
For otherwise the short roots
orthogonal to $q$ project into some triple $T$ of mutually
orthogonal antipodal pairs of short vectors of ${V}$. If $q$
lies on every {cross} in $\calC$ then there is a way to choose a
root in $q^\bot$ for each $C\in\calC$, such that the image in
${V}$ of
the root is one of the vectors of $C$. In particular, there is a
way to choose an element $v\in C$ for each $C\in\calC$, such
that
(1) if $v$ is short then $v\in T$, and (2) the span of
all the $v$'s has dimension at most $4$. For each of the 540
possibilities for $T$ one can count the number of ways to
choose vectors $v$ satisfying (1) and (2). It turns out that there
are no ways to make such a choice, and it follows that $q$
cannot lie on only 3 (or fewer) short mirrors.
\smallskip
We have shown that the short roots of $q^\bot$ project onto some
quadruple of mutually orthogonal antipodal pairs of short
vectors of ${V}$, which we will denote by $T$. As in
the previous paragraph, there is a way to choose an element
$v\in C$ for each $C\in\calC$, such that (1) and (2) are
satisfied. For each of the 135 possibilities for $T$, one can
count the number of ways to make such a set of choices. It turns
out that for only one quadruple is there a way to do this, 
and this quadruple consists of
$\pm(0,1,0,0,0),\ldots,\pm(0,0,0,0,1)$. Therefore the images in ${V}$
of the short roots of $q^\bot$ are these 8 vectors. Since these
are also the images of the short roots of $p^\bot$, the
$\Gamma$-equivalence of $p$ and $q$ follows from Lemma~\foo.
\smallskip
Now we turn to the boundary points of $X$. If $b$ is a boundary
point of $\calB^*_4$ then we may represent it by a primitive
isotropic lattice vector $w$, and a cross contains $b$ just if
it contains the image in $V$ of a root orthogonal to $w$. One
can check that every nonisotropic element of $V$ that is
orthogonal to the image $\bar w$ of
$w$ is the image of a root in $w^\perp$. (This is
easy to check for any given $w$, and the result follows for
general $w$ because of the transitivity of $\Aut\Lambda$.) It
follows that the set of crosses $\calC_b$ containing $b$
consists of the crosses which contain a vector of $V$ orthogonal
to $\bar w$. It is easy to compute the sets $\calC_b$ for each
of the~40 orbits of boundary points, and to check that no
$\calC_b$ is a subset of $\calC$. This proves the second part of
the theorem. The first part then follows, because no point of $X$
except for the image of $p$ lies on every cross in $\calC$, and
this point lies on only~69 of the ~135 crosses.
\smallskip
Now we show that no point of $\calB/\Gamma$ lies on all the
crosses in $\calC_b$, for any boundary point $b$. The proof is
almost identical to the one used above. By symmetry it suffices
to treat just one $\calC_b$. If $q\in\calB_4$, then the short
roots of $q^\perp$ project into some quadrouple $T$ of mutually
orthogonal short vectors of $V$. If $q$ lies on every cross in
$\calC_b$ then there is a way to choose an element $v\in C$ for
each $C\in\calC_b$, such that (1) and (2) are satisfied.
An enumeration shows that there is no way to make such a choice,
and the claim follows. Finally, it is easy to compare the
$\calC_b$'s with each other as $b$ varies over the boundary
points, and check that none of the $\calC_b$'s contains any
other. It follows that for each boundary point $b$ of $X$, $b$ is the
only point of $X$ that lies on all the crosses containing $b$.
This completes the proof.
\smallskip
We verified the enumerations with a computer program written in C++,
which ran to completion
in less than a minute. Repeatedly checking condition (2) required more
than $4\times10^8$ row-reduction operations, and we did this
efficiently by enumerating the $3^5$ elements of ${V}$ and
preparing a lookup table of all $3^{5\cdot2}$ possible
row-reductions. \qed
\neupara{Automorphic forms on the ball}%
Borcherds has given two constructions for automorphic forms on
$\O(2,n)$, which we will use to build automorphic forms on the
4-ball. Here we will use his additive lift [Bo1,\S14], which
generalizes correspondences of Shimura, Doi-Naganuma, Maass,
Gritsenko, and others. In the next section we will discuss his
other construction, which uses infinite products.
\smallskip
We begin in the setting of section~3, 
with $\goto$ an order in
an imaginary quadratic number field, $L$ an $\goto$-lattice of
signature $(1,n)$, $\calB(L)$ the associated ball in projective
space, and $\calB^*(L)$ the union of the ball with the cusps. We
assume that $L$ is integral (all inner products lie in $\goto$)
and that $n>1$, so that $L$ has dimension at least~3. 
We define $\tilde\calB(L)$ and $\tilde\calB^*(L)$ to
be the preimages of $\calB(L)$ and $\calB^*(L)$ in
$L\tensor_{\goto}\cz$. If $G$ is a subgroup of 
$\Aut(L)$ and $v:G\to S^1\subset\cz^\bullet=\cz-\{0\}$ is a character of
$G$ then an automorphic form of weight $k\in\gz$ with respect
to $G$ and $v$ is a holomorphic function
$f:\tilde\calB(L)\to\cz$ satisfying
\medni
\item{\rm a)} $f(tz)=t^{-k}f(z)$ for $t\in\cz^\bullet$, and 
\item{\rm b)} $f(\gamma z)=v(\gamma)f(z)$ for $\gamma\in G$.
\medni
(If $n$ were~1, so that $\calB(L)$ were one-dimensional, then we
would impose an additional condition of regularity at the cusps.)
We denote the space of all such forms by $[G,k,v]$, or by $[G,k]$
if $v$ is trivial.
\smallskip
One can extend an automorphic form $f:\tilde\calB(L)\to\cz$ to
$\tilde\calB^*(L)$ in a natural way, providing boundary
values for $f$. If $a$ is an isotropic element of
$\tilde\calB^*(L)$, so that it represents a cusp, then by the
non-degeneracy of $\spitz{\cdot,\cdot}$ we may choose $b\in
L\tensor_{\goto}\cz$ satisfying $\spitz{a,b}\neq0$. For all
$\tau\in\cz$ with sufficiently large imaginary part, $\tau
a+2\imag\spitz{a,b}b$ has positive norm. The limit
$$f(a):=\lim_{\Im\tau\to\infty}f(\tau a+2\imag\spitz{b,a}b)$$
exists and is independent of the choice of $b$.
This follows from the Fourier Jacobi expansion
of $f$ at a cusp; we refer to [Sh] for more
details.
\smallskip
An automorphic form $f\in[G,k,v]$ is of course not a function on
$\calB(L)$ unless  $k=0$. But it is clear that the zero-locus of
$f$ is preserved by $G$ and scalar multiplication, so the
set of zeros of $f$ in $\calB^*(L)/G$ is well-defined. It is a
closed algebraic subvariety of pure codimension one.
\zwischen{Borcherds' additive lift}%
We consider the $\gz$-lattice $M$ underlying $L$, which is the
underlying $\gz$-module equipped with the even integral
bilinear form
$$(a,b):=\spitz{a,b}+\spitz{b,a},$$
which has signature $(2,2n)$. The dual lattice with respect to
$(\cdot,\cdot)$ is denoted $M'$, and $M'/M$ is a finite
group. We remark that if $\alpha,\beta\in M'/M$ then $(\alpha,\alpha)$
and $(\beta,\beta)$ are well-defined modulo~2, while 
$(\alpha,\beta)$ is well-defined modulo~1.
The group $\SL(2,\gz)$ acts on the group ring
$\cz[M'/M]$ by means of the Weil representation $\varrho_M$,
which is defined in terms of the standard generators
$$T=\pmatrix{1&1\cr0&1},\quad S=\pmatrix{0&-1\cr1&0}$$
by
$$\eqalign{\varrho_M(T)e_\alpha&=\exp(\pii(\alpha,\alpha))e_\alpha,\cr
  \varrho_M(S)e_\alpha&={\imag^{n-1}\over\sqrt{|M'/M|}}
  \sum_{\beta\in M'/M}\exp(-2\pii(\alpha,\beta))e_\beta.}$$
(We denote the standard generators of the group ring $\cz[M'/M]$
by $e_\alpha$, with $\alpha$ varying over $M'/M$.)
The Weil representation factors through $\SL(2,\gz/N\gz)$,
where $N$ is the smallest natural number such that
${N\over2}(a,a)$ is integral for all $a\in M'$.
\smallskip
The inputs of Borcherds' additive lift are vector valued modular forms
$f:H\to\cz[M'/M]$
on the usual upper half plane $H$
with respect to the Weil representation.
More precisely, we require that $f=(f_\alpha)_{\alpha\in M'/M}$ 
satisfy
$$\leqalignno{
&f_\alpha(\tau+1)=e^{\pii(\alpha,\alpha)}f_\alpha(\tau),&1.\cr
&f_\alpha\Bigl(-{1\over\tau}\Bigr)=\tau^{k+1-n}
	{\imag^{n-1}\over\sqrt{|M'/M|}}
     \sum_{\beta\in M'/M}
e^{-2\pii(\alpha,\beta)}
      f_\beta(\tau)\hbox{, and}
	&2.\cr
	&\hbox{\rm $f$ is holomorphic at the cusp
	infinity.\qquad}
	&3.}$$
Borcherds' additive lift allows also inputs which have
poles at the cusps, but we do not need this extension. 
But even in the case of  modular forms
which are regular at the cusps, Borcherds extended
previous constructions because he imposes no restriction on the weight
of $f$, and does not require that $f$ be a cusp form.
\smallskip
The additive lift is a linear map $\Psi$ from the space of such
$f$ into a certain space of automorphic forms on $\calB(L)$. We give
its important properties  in the following
theorem, which is a specialization of
Theorem 14.3 in [Bo1] to $\U(1,n)\subseteq\O(2,2n)$.
\proclaim
{Theorem}
{Let $G$ be the subgroup of $\Aut(L)$ that acts
trivially on $M'/M$.
There exists a linear map $\Psi$ (the additive lift) from the space
of elliptic modular forms with the properties 1--3 above
into the space $[G,k]$ of automorphic forms of weight $k$
with respect to $G$ and the trivial character.
This lifting is equivariant with respect to
the action of $\Aut(L)$.
($\Aut(L)$ acts on $[G,k]$ because
$G$ is normal in $\Aut(L)$, and on the space
of elliptic modular forms via its action on $M'/M$.)
}
Badd%
\finishproclaim
Furthermore, Borcherds shows how to compute the 
values of $\Psi(f)$ at the cusps of $\tilde\calB^*(L)$
from the
Fourier coefficients of $f$.

\smallskip
We now turn to the case of interest, with $\goto=\calE$ and
$L=\Lambda$. The $\gz$-lattice underlying
the 1-dimensional lattice $\calE$ is the $A_2$ root lattice (the
hexagonal lattice in the plane with minimal norm~2), which has
index~3 in its dual. From the definition of $\Lambda$ as a
direct sum, we see that $M'/M$ has order $3^5$. Indeed more is
true: $M'$ coincides with $(\sqrt{-3})^{-1}\Lambda$, so that
$M'/M$ is canonically isomorphic to the $\fz_3$-vector space
$V=\Lambda/\sqrt{-3}\,\Lambda$ introduced in section~2. 
In
particular, $G$ is the congruence subgroup $\Gamma$. One can
check that if $\alpha,\beta\in M'/M$ then $(\alpha,\beta)$ is
$0$, $2/3$ or $-2/3$ (modulo~1) according to whether the
corresponding elements of $V$ have inner product $0$, $1$ or
$-1$ (in $\fz_3$). Similarly, if $\alpha\in M'/M$ then
$(\alpha,\alpha)$ is $0$, $2/3$ or $-2/3$ (modulo~2) according
to whether the corresponding element of $V$ has norm $0$, $1$ or
$-1$. It follows from this that the level of the Weil
representation is $N=3$, so that the representation factors through
$\SL(2,\fz_3)$. We will usually write $V$ in place of $M'/M$
to lighten the notation.
\smallskip
We apply \Badd\ in the simplest case, where $f$ is a modular form
of weight $0$, hence a constant, which is to say an element of
$\cz[V]^{\SL(2,\fz_3)}$.
The weight being $0$ means that the exponent $1-k+n$ of $\tau$ in the
second transformation rule is $0$, so that $k=n-1=3$. Therefore
Borcherds' additive lift gives a linear map
$$\cz[V]^{\SL(2,\fz_3)}\lo [\Gamma,3].$$
We remark that since $\Gamma$ contains the cube roots of unity
acting as scalars, every automorphic form on $\tilde\calB_4$ for
$\Gamma$, with trivial character has weight divisible by~3.
Our first task is to find some elements of
$\cz[V]^{\SL(2,\fz_3)}$.
\proclaim
{Lemma}
{Let $a_0,\ldots,a_4$ be an orthogonal basis for $V$ consisting
of one long vector and four short vectors, and let
$C=(C_\alpha)_{\alpha\in V}\in\cz[V]$ be defined by the
condition that $C_\alpha$ is the complex number $1$, $0$ or $-1$ according to
whether $\prod_i(\alpha,a_i)$ is the element $1$, $0$ or $-1$ of
$\fz_3$. Then $C$ is invariant under the Weil
representation. Furthermore, $C$ changes sign under reflection
in any of the $a_i$, and is characterized up to a scalar by
this property.
}
constr%
\finishproclaim
To avoid the impression that the $C$'s were
discovered by clever guesswork, we should mention that we found this
construction quite late, following extensive computer work.
\smallskip
{\it Proof.}
The behavior of $C$ under the reflections is obvious, and the
invariance under $\SL(2,\fz_3)$ may be checked by a computer
calculation. To see the last claim, suppose
$D=(D_\alpha)\in\cz[V]$ has the stated property. If $\alpha$
is orthogonal to one of the $a_i$ then we have
$D_\alpha=-D_\alpha$ by the transformation rule, so that
$D_\alpha=0$. The remaining $\alpha$ fall into a single orbit of
size~32 under the group $(\gz/2)^5$ generated by the reflections
in the $a_i$, so all the remaining $D_\alpha$ are determined by
any one of them.
\qed
\smallskip
It is easy to work this out explicitly in an example: if
$a_0,\ldots,a_4$ are $(1,0,\ldots,0),\ldots,(0,\ldots,0,1)$ then
$C$ is supported on those $\alpha$ of the form
$(\pm1,\ldots,\pm1)$, with $C_\alpha=+1$ or $-1$ according to
whether there are an even or odd number of minus signs. Note
that $C$ is supported on the isotropic vectors in $V$, which is
not immediately obvious from the construction. 
It
follows from the lemma and Theorem~\Badd\ that to each cross there
is associated an automorphic form on $\calB_4$,
well-defined up to sign. We will see below that the zero-locus
of this form is exactly the associated cross in $\calB_4$. To
resolve the sign ambiguity it is convenient to introduce the
notion of a {\it signed cross}. This is just a basis
$\{a_0,\ldots,a_4\}$ as in the lemma, modulo the equivalence
relation that $\{a_0,\ldots,a_4\}\sim\{a_0',\ldots,a_4'\}$ if
the $a_i'$ differ from the $a_i$ by a permutation and evenly
many sign changes. It is clear that there are two signed crosses
for every cross, and that the lemma assigns an element of
$\cz[V]$ to each of the~270 signed crosses.
\proclaim
{Lemma}
{The space 
$$\cz[M'/M]^{\SL(2,\gz)}=\cz[V]^{\SL(2,\fz_3)}$$
has dimension $10$ and is spanned by the elements of $\cz[V]$
associated to the signed crosses. The group $\O(5,3)$ acts
irreducibly on this space, with $W(E_6)$ acting by its unique
10-dimensional irreducible representation and the central
involution acting by $-1$. The multiplicity of this
representation in $\cz[V]$ is one.
}
Sing%
\finishproclaim
{\it Proof.} 
It is easy to make a computer construct the elements $C$ of
$\cz[V]$ associated to the signed crosses and check that their
complex span $Z$ is 10-dimensional. Consulting the character
table shows that any 10-dimensional representation of $W(E_6)$
is either trivial, or the irreducible representation in 10
dimensions, or else the sum of the (unique) irreducible
6-dimensional representation and a 4-dimensional trivial
one. These may be distinguished by the trace of almost any group
element, say a short reflection $R$, which has ATLAS [C] conjugacy
class 2C. The fixed space of $R$ in $Z$ is spanned by the
vectors $C+R(C)$ where $C$ is as above. It is easy to check that
this space has dimension 5, so that $R$ has trace 0, so that
$W(E_6)$ acts irreducibly on $Z$. It is obvious that each $C$
changes sign under the central involution of $\O(5,3)$.
If the multiplicity of this $\O(5,3)$-representation in $\cz[V]$
were more than one, then the subspace of $\cz[V]$ that changed
sign under the reflections of each vector in a cross
would have dimension~$>1$, contrary to Lemma~\constr.
\smallskip
To see that $Z$ is all of $\cz[V]^{\SL(2,\fz_3)}$, suppose
$C=(C_a)_{a\in V}$ is an element of
$\cz[V]^{\SL(2,\fz_3)}$.  Invariance under $T$ means
that $C$ is supported on the $81$ isotropic elements.
Invariance under $S^2=-E$ means $C_{-a}=-C_a$.
Invariance under $S$ can be read as a linear equation in $40$
indeterminates, and it is easy to make a computer check that the
space of solutions has only~10 dimensions.
If one is prepared to do more work with group characters, one
can of course find the complete decomposition of $\cz[V]$ under
$\SL(2,\fz_3)\times \O(5,3)$; this is done in [Fr].
\qed
\smallskip
We will write $W$ for the image of Borcherds' additive lift
$\cz[V]^{\SL(2,\fz_3)}\to[\Gamma,3]$. 
Our next theorem asserts that the automorphic forms
we have constructed are nontrivial:
\proclaim
{Proposition}
{Borcherds' additive lift 
$$\cz[V]^{\SL(2,\fz_3)}\to W\subset[\Gamma,3]$$
is an $\O(5,3)$-equivariant embedding.
}
addEm%
\finishproclaim
{\it Proof.}
The $\O(5,3)$-equivariance is part of Theorem~\Badd. To prove
injectivity, we construct an inverse by using the boundary
values of the automorphic forms. Namely, if
$f\in[\Gamma,3]$ then we define $C=(C_\alpha)_{\alpha\in V}$ by
taking $C_\alpha=0$ if $\alpha$ is zero or nonisotropic, and
$C_\alpha=f(\tilde\alpha)$  otherwise, where $\tilde\alpha$ is
any primitive isotropic vector in $\Lambda$ representing
$\alpha$. This definition is independent of the choice of $\tilde\alpha$
because $f$ is $\Gamma$-invariant and all the primitive
isotropic preimages of $\alpha$ are $\Gamma$-equivalent
(Lemma~\Orb). The irreducibility of $\cz[V]^{\SL(2,\fz_3)}$ as an
$\O(5,3)$-module and the fact that its multiplicity in $\cz[V]$
is one imply that the composition
$$
\cz[V]^{\SL(2,\fz_3)}
\;
\buildrel \hbox{\rm \sevenrm additive lift} 
\over{\hbox to 55pt{\rightarrowfill}}
\;
W
\subset
[\Gamma,3]
\;
\buildrel\hbox{\rm \sevenrm boundary values} 
\over{\hbox to 75pt{\rightarrowfill}}
\;
\cz[V]
$$
is a scalar. The problem is to show that this scalar is nonzero.
This is a straightforward but tedious calculation using
Borcherds' formulae for the Fourier expansions of
additive lifts ([Bo1], 14.3) and the explicit embedding
of $\calB_4$ into the hermitian symmetric
space of $\O(2,8)$. The latter space  consists of two-dimensional
positive definite real subspaces of $M\tensor_\gz\rz$. Every positive definite
complex line in $\Lambda\tensor_{\calE}\cz$ 
(i.e.\ a point in $\calB_4$) defines
such a subspace. One has to express this embedding
in the coordinates which Borcherds uses in his theorem 14.3. 
Details of this calculation can be found in section~6 of [Fr].
\qed
\smallskip
Lemma~\Sing\ shows that our~270 automorphic forms satisfy many
linear equations. Some of these are easy to see, and those
treated in the following lemma will receive an elegant
geometric interpretation in section~7. 
To formulate the lemma we note that there is an
$\O(5,3)$-invariant inner product on $W$, which is
unique up to scale, by the irreducibility of the representation.
\proclaim
{Lemma}
{Let $v$ be a long vector of $V$. Then the automorphic 
forms associated to the six signed crosses involving $v$ lie in
a 2-dimensional subspace of $W$, and form a scaled copy of the $A_2$
root system, i.e., the vertices of a regular hexagon centered at
$0$.
}
EinF%
\finishproclaim
{\it Proof.}
One can check this by computing the inner products of the 6
elements of $\cz[V]$, using the restriction of the
inner product 
$$
\bigl((C_\alpha),(D_\alpha)\bigr)=\sum_{\alpha\in V}C_\alpha D_\alpha,
$$
which is obviously $\O(5,3)$-invariant and therefore the natural
inner product.  But here is a better argument.  The reflection $R$
in $v$ is not in the simple subgroup of $W(E_6)$, but $-R$ is,
and has conjugacy class $2A$ in ATLAS notation. Consulting the
character table shows that $-R$ has trace $-6$, so that the
subspace $Z$ of $W$ that $R$ negates has
dimension~2.
Lemma~\constr\ associates to each of
the three crosses a one-dimensional subspace of $Z$, with
two generators coming from the associated signed crosses. The
reflection of $V$ in any short root of one of these crosses preserves
that cross, acts as $-1$ on the the associated subspace, and
exchanges the other two crosses and hence the corresponding
subspaces. 
Therefore the three subspaces meet each other at angles of
$\pi/3$, the subgroup of
$\O(5,3)$ generated by the reflections in the short vectors
orthogonal to $v$ acts on $Z$ by the $A_2$ Weyl group, and
the~6 elements of $Z$ coming from the crosses form a
copy of the $A_2$ root system.
\qed
\smallskip
We recall the notion of the divisor of an automorphic form.
Let $Y\subset X$ be
an irreducible subvariety of codimension one.
We denote by $e_Y$ the ramification degree with respect to
the natural projection $\pi:\calB_4\to X$
(counted as $1$ if $\pi$ is unramified along $Y$). If $Y$ is a short
mirror this ramification degree is three, 
because the triflections are contained in $\Gamma$.
For any other $Y$, such as a long mirror, it is one.
If $F$
is a nonzero automorphic form for $\Gamma$ and the trivial
character, 
then the vanishing order of $F$ along $\pi^{-1}(Y)$
is divisible by $e_Y$. We call the quotient of
this vanishing order by $e_Y$ the order of $F$ along $Y$ and
denote it by $n_Y(F)$. The divisor of $F$ in $X$ is the finite sum
$$(F):=\sum_{Y\subset X}n_Y(F)Y.$$
We consider a {cross} in $X$ as a divisor with multiplicity one
at all 5 of its mirrors.
A fundamental result for this paper is
\proclaim
{Theorem}
{Let $F\neq0$ lie in the one dimensional space of automorphic forms
associated to a cross. Then the divisor of $F$ in $X$ is exactly
this cross. The 270 automorphic forms associated to the signed
crosses have no common zeros in $\calB^*_4$.}
farD%
\finishproclaim
\smallni
To prove this we will need a result whose proof we postpone
to the next section. We remark that the form $\chi_4$
given here was first discovered by Borcherds [Bo3].
\proclaim
{Theorem}
{There are automorphic forms $\chi_4\in[\Aut(\Lambda),4,v]$ and
$\chi_{75}\in[\Aut(\Lambda),75,v']$, for some characters $v$ and $v'$ of
$\Lambda$, such that the divisors of $\chi_4$ and $\chi_{75}$ in
$\calB_4$ are the sum of the short mirrors and the sum of the
long mirrors, respectively, with multiplicity one.
}
farDI%
\finishproclaim
{\it Proof of Theorem~\farD.}
Suppose the cross is $\{\pm a_0,\ldots,\pm a_4\}\subset V$. If
$\tilde a$ is a root of $\Lambda$ representing any of the $\pm
a_i$, and $R$ is the biflection in $\tilde a$, then the relation
$F\circ R=-F$ (which follows from the construction of $F$)
implies that $F$ vanishes along the mirrors of $\tilde
a$. Furthermore, if $\tilde a$ is a short root then $F$ is
invariant under the triflection in $\tilde a$, so that the
multiplicity in $\calB_4$ is at divisible by~3. It
follows that the divisor of $F$ in $\calB_4$ contains the
short mirrors of the cross with multiplicity~3, plus the
long mirrors of the cross. To prove the theorem it suffices to
show that this is the full divisor of $F$. To see this we
construct the product $P$ of all~270 automorphic forms, and divide
$P$ by $\chi_4^{90}\chi_{75}^6$, where $\chi_4$ and $\chi_{74}$
are as in Theorem~\farDI. The
quotient is holomorphic because $P$ vanishes to order at least
$6$ along each long mirror in $\calB_4$ and least
$270\cdot3\cdot4/36=90$ along each short mirror. The quotient
has weight $270\cdot3-90\cdot4-75\cdot6=0$, so is constant. It
is nonzero because each $F$ is nonzero. Therefore
the divisor of $P$ is the same as that of
$\chi_4^{90}\chi_{75}^6$; since this is also the sum of the
``known'' divisors of the various $F$, the first statement of
the
theorem follows. The second follows immediately from this and
Theorem~\Fraemp.
\qed
\neupara{Borcherds products (and proof of theorem \farDI)}%
We recall some facts about
automorphic forms on $\O(2,n)$, where $\O(2,n)$ is the
orthogonal group of a real vector space $V$ with a symmetric
bilinear form $(\cdot,\cdot)$ of signature $(2,n)$.
Let $\calH_n$
denote the hermitian symmetric space associated to $\O(2,n)$.
It can be realized as an open subset of the quadric
defined by $(z,z)=0$ in the projective space $P(V(\cz))$,
where we extend $(\cdot,\cdot)$ to a $\cz$-bilinear form
on $V(\cz)$. Namely, it is one of the two connected
components of the open subset defined
by $(z,\bar z)>0$. A subgroup $\O'(V)$ of index two of the
orthogonal group $\O(V)$ acts biholomorphically on $\calH_n$. 
Let $\tilde\calH_n$ denote the
inverse image of $\calH_n$ in $V(\cz)$. 
We restrict henceforth to the case $n>2$ for convenience.
If $M$ is an even integral $\gz$-lattice in $V$, then a  
meromorphic automorphic
form of weight $k\in\gz$ with respect to a subgroup $G$ 
of finite index in $$\O'(M)=\O(M)\cap\O'(V)$$
and a character $v$ of $G$ is a meromorphic
function $f$ on $\tilde\calH_n$ with the properties
\smallni
a) $f(\gamma z)=v(\gamma)f(z)$ for all $\gamma\in G$.
\hfill\break
b) $f(tz)=t^{-k}f(z)$ for all $t\in\cz^\bullet$.
\smallni
\smallskip
We next recall the notion of a Heegner divisor: let $m$ be a
negative rational number and let $\alpha$ be an element of
$M'/M$, where $M'$ denotes the dual lattice. The Heegner divisor
$H(\alpha,m)\subset\calH_n$ is the union of the orthogonal
complements $v^\perp\cap\calH_n$ where $v$ runs through all
elements of $M'$ satisfying
$$v\equiv\alpha\;\mod\; M\quad\hbox{\rm and}\quad (v,v)=2m.$$
We consider $H(\alpha,m)$ as a divisor on $\calH_n$ by attaching
multiplicity $1$ to all components. It is obvious that 
$H(\alpha,m)=H(-\alpha,m)$, so that the divisor depends
only on $m$ and the image of $\alpha$ in $(M'/M)/\pm1$.
\smallskip
Borcherds introduced in [Bo1] a method for constructing
automorphic forms on $\calH_n$ whose divisors are sums of
Heegner divisors. Then, in
[Bo2], he constructed a `space of obstructions'
to the use of this technique for constructing automorphic forms
with divisor equal to some given sum of Heegner divisors. This space
consists of all elliptic modular forms of weight 
$$
k:=(2+n)/2
$$
with respect to the dual $\varrho^*:=\varrho_M^*$ of the Weil
representation.  We restrict to the case of even $n$ since the
Weil representation simplifies and this is the only case we
need.  Such a form $(f_\alpha)_{\alpha\in M'/M}$ is required to
be holomorphic at the cusp at infinity and satisfy the
transformation laws
$$\leqalignno{
&f_\alpha(\tau+1)=e^{-\pii (\alpha,\alpha)}f_\alpha(\tau)&1.\cr
&f_\alpha\Bigl(-{1\over\tau}\Bigr)=-\sqrt{\tau\over\mag}^{2+n}
	{1\over\sqrt{|M'/M|}}\sum_{\beta\in M'/M}e^{2\pii(\alpha,\beta)}
      f_\beta(\tau).
	&2.}
$$ 
As in section~4, 
we note that $(\alpha,\alpha)$ is well-defined
modulo~2 and $(\alpha,\beta)$ is well-defined modulo~1, so that
these formulas make sense.
\smallskip
Elements of the space of obstructions can be constructed
by means of Eisenstein series, as follows. We write $R$ for the group ring
$\cz[M'/M]$ and $R_0$ for the subspace on which
$(-1)^k\varrho^*(-E)$ acts trivially. Since $(-1)^k\varrho^*(-E)$
acts by exchanging $e_\alpha$ and $e_{-\alpha}$, where the
$e_\alpha$ form the standard basis of $R$ as $\alpha$ varies over
$M'/M$, a basis for $R_0$ is given by the elements 
$$
e_\alpha+e_{-\alpha},\qquad\alpha\in(M'/M)/\pm1.
$$
If $\xi\in R_0$ satisfies 
$$\varrho^*\pmatrix{1&1\cr0&1}\xi=\xi$$
then $(c\tau+d)^{-k}\varrho^*(Q)^{-1}\xi$ remains unchanged if one
replaces $Q$ by $PQ$, where $P=\pm\smallmatrix{1}{n}{0}{1}$ lies in
the stabilizer $\SL(2,\gz)_\infty$ of $\infty$ and
$Q=\smallmatrix{a}{b}{c}{d}\in\SL(2,\gz)$. This lets us define the
Eisenstein series 
$$
E_\xi(\tau)=\sum_{Q=\smallmatrix{a}{b}{c}{d}\in\SL(2,\gz)_\infty\backslash\SL(2,\gz)}
(c\tau+d)^{-k}\varrho^*(Q)^{-1}\xi\;,
$$
which is a modular form of weight $k$ with respect to
$\varrho^*$, so it lies in the space of obstructions. 
Furthermore,
$$\lim_{\im\tau\to\infty} E_\xi(\tau)=\xi.$$
In particular, if $\xi\neq0$ then $E_\xi$ is not a cusp form.
\proclaim
{Remark} 
{Under our assumption $n>2$ we have $k>2$, and in this
case the space of Eisenstein series of weight $k$ and with
respect to $\varrho^*$ is isomorphic the space of all $\xi\in V_0$
with
$$\varrho^*\pmatrix{1&1\cr0&1}\xi=\xi\qquad\left(\schreib{and }
\varrho^*(-E)\xi=(-1)^k\xi\right)\>.$$
}
dimeis%
\finishproclaim
Using Remark~\dimeis\ 
one can reformulate a fundamental result of
Borcherds [Bo1], [Bo2] as follows.
\proclaim
{Theorem}
{Suppose $D$ is a finite $\gz$-linear combination
$$\sum_{\alpha\in (M'/M)/\pm1,\;m<0}  C(\alpha,m)H(\alpha,m)$$
of Heegner divisors and
$G$ is the subgroup of $\O'(M)$ that acts trivially on
$M'/M$. Then $D$
is the divisor of a meromorphic automorphic form on $\calH_n$
with respect to some character of $G$ if for every cusp
form $f$ in the space of obstructions, say
$$
f_\alpha(\tau)=\sum_{m\in\qz} a_\alpha(m)\exp(2\pii m\tau),
\eqno(5.1)
$$
the relation
$$
\sum_{m<0,\,\alpha\in M'/M}a_\alpha(-m)C(\alpha,m)=0
\eqno(5.2)
$$
holds.
The weight of such an automorphic form is
$$ \sum_{m\in\qz,\,\alpha\in M'/M}b_\alpha(m)C(\alpha,-m),$$ 
where $b_\alpha(m)$ denotes the Fourier coefficients
of the (unique) Eisenstein series with constant term
$$b_\alpha(0)=\cases{-1/2&if $\alpha=0$,\cr0&otherwise.}$$
}
bopro%
\finishproclaim
\smallskip
We want to apply this theorem to our lattice 
$\Lambda=\Eisenstein^{1,4}$, or rather to its underlying $\gz$-lattice $M$.
The obstructions have weight $k=(2+8)/2=5$, and if the space of
obstructions vanished then the existence of the forms of
Theorem~\farDI\ would follow easily from Theorem~\bopro\ by restriction
from $\calH_8$ to $\calB_4$. There are obstructions, and
even cuspidal obstructions, but we will show that
there are no $\O'(M)$-invariant cusp forms in the space of
obstructions, and this turns out to be enough to establish
Theorem~\farDI. Here, $\O'(M)$ acts via its action on $M'/M$.
\proclaim
{Theorem}
{For $M$ equal to the $\Z$-lattice underlying
$\Lambda=\Eisenstein^{1,4}$, the space of $\O'(M)$-invariant
obstructions has dimension two and is spanned by Eisenstein
series. The space of invariant cuspidal obstructions vanishes.}
noobstr%
\finishproclaim
{\it Proof.}
The $\O'(M)$-invariant part of $\cz[M'/M]$ has dimension~4,
because $\O'(M)$ acts with~4 orbits (or `types') on
$M'/M$. 
The {\it type\/} of an element $\alpha\in M'/M$ is defined
as $00$ if $\alpha$ is the zero element and as $t\in\{0,1,2\}$
if $\alpha$ is different from zero and 
$(\alpha,\alpha)\equiv 2t/3$ mod $2$. There are~1, 80, 90 and~72
elements of $M'/M$ of types~00, 0, 1 and~2, respectively. We will
express an invariant obstruction $h$ as $(h_{00},h_0,h_1,h_2)$,
where each $h_t$ is the sum of the $h_\alpha$ as $\alpha$ varies
over the elements of $M'/M$ of type $t$. A calculation allows
one to determine the action of $\SL(2,\gz)$ with respect to this
basis. It turns out that the standard generators
$T=\smallmatrix1101$ and $S=\smallmatrix0{-1}10$ act by 
$$
\varrho^*(T)=\pmatrix{1&&&\cr&1&&\cr&&\omega^2&\cr&&&\omega}
\quad\hbox{and}\quad
\varrho^*(S)={\imag\over3^{5/2}}\pmatrix{	1&1&1&1\cr
					80&-1&8&-10\cr
					90&9&-9&0\cr
					72&-9&0&9\cr}
\;.
$$
Borcherds [Bo2] gives a formula for the dimension of the space
of elliptic modular forms of given weight that tranform
according to some given representation of $\SL(2,\gz)$, in terms
of the eigenvalues of certain elements of $\SL(2,\gz)$. Applying
this formula to the 4-dimensional representation above shows
that the space of obstructions has dimension~2. On the other
hand, the space of Eisenstein series is also 2-dimensional,
because another calculation shows that the subspace of
$\cz[M'/M]^{\O'(M)}$ whose elements satisfy the conditions of
Remark~\dimeis\ is 2-dimensional. Since a cuspidal Eisenstein
series vanishes identically, the theorem follows.
\qed
\proclaim
{Corollary} {With $M$ as in Theorem~\noobstr, every divisor $D$
as in Theorem~\bopro\ that is $\O'(M)$-invariant is the divisor
of a form on $\calH_8$ that is automorphic with respect to some
character of $\O'(M)$.}  
formsexist%
\finishproclaim
{\it Proof.}
Since $D$ is $\O'(M)$-invariant it satisfies condition (5.2) 
for all $f$ as in (5.1) 
if and only if it satisfies the condition
for all such $f$ that are also $\O'(M)$-invariant. Therefore the theorem
assures us of the existence of an automorphic form for
$G\subseteq\O'(M)$ with divisor $D$, and since $D$ is
$\O'(M)$-invariant the form must be automorphic with respect to
some character of $\O'(M)$ itself.
\qed
\smallskip
In order to find the weights of the forms constructed in this
way, we need the Eisenstein series with constant term
$b_\alpha(0)=-1/2$ for $\alpha=0$ and $b_\alpha(0)=0$ for other
$\alpha$. To compute this series we construct a basis for the
space of Eisenstein series and then take a suitable linear
combination. 
The Weil representation
factors through $\SL(2,\gz/3\gz)$, so
our Eisenstein series are linear combinations of the classical
elliptic Eisenstein series of level $3$, namely
$$G_k(\tau;c,d,N):=\mathop{{\sum}'}_{{m\equiv c\atop n\equiv d}\mod N}
 {1\over(m\tau+n)^k}\;,$$
where the level $N$ is~3 and the weight $k$ is~5.
We write
$E_1$, $E_2$, $E_3$ and $E_4$ for the four classical Eisenstein series
corresponding to the values
$(c,d)=(0,1)$, $(1,0)$, $(1,1)$ and $(1,2)$.
We will continue to use the notation introduced in the proof of
Theorem~\noobstr. 
\proclaim
{Proposition}
{With $M$ as in Theorem~\noobstr, a basis for the space of
$\O'(M)$-invariant obstructions consists of the Eisenstein
series $f$ and $g$ given by
$$
\eqalign{f_{00}&={\imag\sqrt3\over18}E_1-{1\over18}(E_2+E_3+E_4)\cr
f_0&=-{5\imag\sqrt3\over 9}E_1+{5\over9}(E_2+E_3+E_4)\cr
f_1&=0\cr
f_2&=E_2+\omega^2E_3+\omega E_4\cr}\quad
\eqalign{g_{00}&={\imag\sqrt3\over18}E_1+{1\over18}(E_2+E_3+E_4)\cr
g_0&={4\imag\sqrt3\over 9}E_1+{4\over9}(E_2+E_3+E_4)\cr
g_1&=E_2+\omega E_3+\omega^2E_4\cr
g_2&=0\;.\cr}$$}
eisseriesbasis%
\finishproclaim
{\it Proof.}
If $h=(h_{00},h_0,h_1,h_2)$ is an $\O'(M)$-invariant Eisenstein
series then each component of $h$ is a $\cz$-linear combination
of $E_1,\ldots,E_4$. The manner in which the $E_i$'s transform
into each other under $\SL(2,\gz)$ is known, and the
transformation laws of $h$ with respect to $\varrho^*$ reduce to a
set of linear conditions on the coefficients of the $E_i$'s. One
simply solves the system of linear equations. (Of course, once
one has the answer one can simply check it.)
\qed
The Fourier coefficients of the Eisenstein series 
can be found in many text books, for example [He, no.~24, section~1] 
or [Fr]. 
This lets one find the Fourier expansions for $f$ and $g$; once
these are known then one can find the unique obstruction $h$
whose Fourier coefficients $b_\alpha(m)$ have constant term
as in Theorem~\bopro. The answer turns out to be given by
$$
\eqalign{
h_{00}&=-1/2+12\,q+225\,{q}^{2}+1092\,{q}^{3}+2892\,{q}^{4}+\cdots\cr
h_0&=1080\,q+16200\,{q}^{2}+87480\,{q}^{3}+260280\,{q}^{4}+673920\,{q}^{5}+\cdots\cr
h_1&=225\,{q}^{2/3}+9360\,{q}^{5/3}+57825\,{q}^{8/3}+219600\,{q}^{11/3}+
540450\,{q}^{14/3}+\cdots\cr
h_2&=12\,{q}^{1/3}+2892\,{q}^{4/3}+28824\,{q}^{7/3}+
112320\,{q}^{10/3}+342744\,{q}^{13/3}+\cdots\cr
}$$
\proclaim
{Proposition}
{With $M$ as in Theorem~\noobstr,
there exists an automorphic form on $\calH_8$ 
for $\O'(M)$, of weight $12$ (resp.\
$225$), whose zeros are the orthogonal complements of the vectors
$v\in M'$ satisfying $(v,v)=-2/3$ (resp.\ $(v,v)=-1/3$).
The vanishing order is one.}
daSi%
\finishproclaim
{\it Proof.}
This follows from Theorems~\bopro and~\noobstr. For the form of weight~12
(resp.~225) we take $D$ to be the sum of the $H(\alpha,m)$ where
$m=-1/3$ (resp. $m=-2/3$) and
$\alpha$ varies over the type~2 (resp.\ type~1) elements of
$(M'/M)/\pm1$.
\qed
\noindent
We note that the form of weight 12 was found
by Borcherds in [Bo3].
\smallskip
{\it Proof of Theorem~\farDI.}
We use the natural embedding  $\U(1,4)\hookrightarrow
\O'(2,8)$ and a compatible holomorphic embedding
$\calB_4\hookrightarrow\calH_8$. The short (resp.\ long) mirrors
in $\calB_4$ are the intersections of the divisors described
in Prop.~\daSi. The orthogonal complements in Prop.~\daSi\ occur in triples
having the same intersection with $\calB_4$, since if $r$ is a
root of $\Lambda$ then $r$, $\omega r$ and $\omega^2r$ are roots with
the same orthogonal complement in $\calB_4$ but different
orthogonal complements in $\calH_8$. Therefore the vanishing
order of the restriction to $\calB_4$ along each short (resp.\
long) mirror is three. Taking
a cube root yields a form of weight $12/3=4$ (resp.\ $225/3=45$).
\qed
\neupara{A model for the moduli space of marked cubic surfaces}%
Recall  the ten dimensional
space $W$ of automorphic forms for $\Gamma$, the congruence subgroup
of level $\sqrt{-3}$ in $\Aut(\Lambda)$.
We know from Theorem~\farD\ that these
forms have no common zero. Therefore, by choosing a basis for
$W$ we obtain an everywhere
holomorphic map
$$\beta:X=\calB_4^*/\Gamma\lo P^9(\cz).$$
This map is algebraic by Chow's theorem. By a result of Hilbert
it is a finite map. Hence the image is a projective
algebraic variety $\calV\subset P^9$ of dimension 4. In fact
more is true:
\proclaim
{Theorem}
{The map 
$\beta:X\to\calV$
is birational.}
BiR%
\finishproclaim
After proving this theorem we will introduce a family of cubic
8-folds, each of which contains $\calV$. Then we will sketch a proof that
these cubic equations actually define $\calV$. Our proof
of~\BiR\ uses only our automorphic forms. In section~7
we will
prove that $\beta$ is actually an embedding, but this relies
heavily on the very extensive calculations and involved
arguments of [Na].
\smallskip
Theorem~\BiR\ follows immediately from the lemma:
\proclaim
{Lemma}
{Let $p$ be the point of $\calB_4$ represented by
$(1,0,0,0,0)\in \Lambda$, and let $\bar{p}$ denote its image in
$X$. Then $\bar p$ is the only point of $X$ mapping to
$\beta(\bar p)$, and $\beta:X\to\calV$ is a local diffeomorphism
at $\bar p$.}
BiRlemma%
\finishproclaim
{\it Proof.}
The first claim is a consequence of the second part of Theorem~\Fraemp. In
order to prove the second claim we will find four elements
of $W$, the sum of whose divisors in $X$ is a normal crossing
divisor at $\bar p$. For this we will need coordinates around
$\bar p$. Coordinates around $p\in\calB_4$ may be taken to be
$z_1,\ldots,z_4\in\cz^4$, with $\sum_i|z_i|^2<1$ as in formula (3.1). 
The 
stabilizer $\Gamma_p$ of $p$, which is
generated (modulo scalars) by the triflections in
the short roots 
$$
\hbox{$(0,1,0,0,0)$, $(0,0,1,0,0)$, $(0,0,0,1,0)$ and
$(0,0,0,0,1)$,}
$$
acts by multiplying the $z_i$ by cube roots of unity.
It follows that local coordinates for $X$ near $\bar p$ are
given by the functions $w_i=z_i^3$. The four short mirrors of
$X$ passing through $\bar p$ are given by the equations
$w_i=0$. The long mirrors in $\calB_4$ that pass through $p$ are
the mirrors of the 216 roots $(0,a_1,\ldots,a_4)$, where two of
the $a_i$ are zero and the others are sixth roots of unity. To
work out their images in $X$ it suffices to treat the case where
the nonzero $a_i$ lie in $\{\pm1\}$, since the orbit of these
under $\Gamma_p$ is the entire set of 216. It is easy to check
that the mirror $z_i=\pm z_j$ in $\calB_4$ projects to the mirror $w_i=\pm
w_j$ in $X$. It follows that in our local coordinates in $X$, the 12
long mirrors through $\bar p$ are given by the equations
$w_i=\pm w_j$ for the various pairs $i\neq j$. 
\smallskip
We claim that for each long mirror $m$ of $X$ passing through
$\bar p$, there is a cross $C$ containing it whose short
mirrors do not pass through $\bar p$. To see this,
consider the three crosses containing $m$. Because only two of
the short mirrors passing through $\bar p$ are orthogonal to
$m$, one of the three crosses contains neither of these
mirrors. Since it cannot contain either of the other short
mirrors, it has the desired property and we take it to be $C$.
Now, it is easy to find four long mirrors $m_1,\ldots,m_4$ 
whose sum is a normal
crossing divisor at $\bar p$, for example those given by
$w_1=\pm w_2$ and $w_3=\pm w_4$. We let $C_i$ be crosses
associated to the $m_i$ as above, and $f_i$ be automorphic forms
associated to the $C_i$. Then the $f_i$ are necessarily linearly
independent, and we extend them to a basis $f_1,\ldots,f_{10}$
of $W$. Of course, one of the $f_i$, say $f_{10}$, does not
vanish at $\bar p$, and then $f_1/f_{10},\ldots,f_9/f_{10}$ are
affine coordinates near $\beta(\bar p)\in P^9$. It follows from
the implicit function theorem and the fact that $f_i$
($i=1,\ldots,4$) has only a simple zero along $m_i$ that $\beta$
is a local diffeomorphism as $\bar p$.
\qed
\smallskip
Next we will find some cubic relations satisfied by our
automorphic forms; these define cubic 8-folds in $P^9(\cz)$ which
contain
$\calV$. It is easy to explain the
origin of these relations: it can happen that there are three
crosses $C_1$, $C_2$ and $C_3$, and another three crosses $C_1'$,
$C_2'$ and $C_3'$, such that as divisors in $X$ they satisfy
$$
C_1+C_2+C_3=C_1'+C_2'+C_3'.
\eqno (6.1)
$$
If $F_i$ and $F_i'$ are nonzero automorphic forms in the
one-dimensional subspaces of $W$ associated to the $C_i$ and
$C_i'$, then the divisors of $F_1F_2F_3$ and $F_1'F_2'F_3'$ are
equal and therefore the two products coincide up to a
scalar. This relation would be trivial if the $C_i'$ were
obtained by permuting the $C_i$, but nontrivial relations do
arise and can be found by studying the geometry of $V$. Here are
some nontrivial cubic relations, which turn out to be the only
ones.  (Only trivial relations can be found if one plays the
same game with pairs rather than triples of crosses.)
\proclaim
{Lemma}
{Let $(a_1,a_2,a_3,b_1,b_2)$ be an ordered 
orthonormal basis of $V$, $S_i$ be
the signed cross given by the basis $\{a_i,a_{i+1}\pm b_1,a_{i-1}\pm b_2\}$,
and $S_i'$ be the signed cross given by 
$\{a_i,a_{i+1}\pm b_2,a_{i-1}\pm b_1\}$, where the subscript of
$a_{i\pm1}$ should be read modulo~3.  Writing $F_i$ and $F_i'$ for the
automorphic forms associated to $S_i$ and $S_i'$, we have 
$F_1F_2F_3=F_1'F_2'F_3'$.}
cubrel
\finishproclaim
{\it Proof.}
We write $C_i$ and $C_i'$ for the
crosses underlying $S_i$ and $S_i'$. 
It is easy to check that (6.1) 
holds, and it follows
that $F_1F_2F_3$ is a constant multiple of $F_1'F_2'F_3'$. 
To determine the constant, 
let $\alpha$ be the isotropic vector $a_1+a_2+a_3\in V$
and let $\tilde\alpha$ be any primitive isotropic element of $\Lambda$
representing $\alpha$. Using the product formula of Lemma~\constr, it is
easy to see that the element of $\cz[V]^{\SL(2,\fz_3)}$ associated to
each $S_i$ and $S_i'$ has component $1$ at $\alpha$. By the
relationship between the values of elements of $W$ at cusps of
$\tilde\calB_4^*$ and at the corresponding elements of $V$ (see the proof
of Prop.~\addEm),
all the $F_i$ and $F_i'$ take the same value at
$\tilde\alpha$. Therefore
$F_1F_2F_3(\tilde\alpha)=F_1'F_2'F_3'(\tilde\alpha)$, and so
$F_1F_2F_3=F_1'F_2'F_3'$. 
\qed
\smallskip
{\it Remarks:}
We will discuss coincidences among these relations, and the fact that
they account for all the relations arising from crosses
$C_i,C_i'$ satisfying (6.1). 
If
$(\hat{a}_1,\hat{a}_2,\hat{a}_3,\hat{b}_1,\hat{b}_2)$ is another
ordered orthonormal basis for $V$, then the relations given by the two bases
are essentially the same if
$$
\{\pm a_1,\pm a_2,\pm a_3\}=\{\pm\hat{a}_1,\pm\hat{a}_2,\pm\hat{a}_3\}
\quad\hbox{and}\quad
\{\pm b_1,\pm b_2\}=\{\pm\hat{b}_1,\pm\hat{b}_2\}.
\eqno (6.2)
$$
By ``essentially the same'' we mean that each relation implies
the other. There are 
$|\O(5,3)|/2^5\,3!\,2!=270$ equivalence classes of ordered orthonormal 
bases under the 
relation (6.2), 
yielding~270 cubic relations.
It is easy to make a
computer enumerate all nontrivial pairs of triples
of crosses $C_i$ and $C_i'$ satisfying (6.1)
and check that every one is a case of our
construction. Therefore we have found all the relations arising
from equalities of sums of triples of crosses.  For convenience
in enumerating the 270 relations, we remark that they are in 1-1
correspondence with the unordered triples of mutually orthogonal
long mirrors in $X$.  To find the relation associated to such a
triple of mirrors, let $a_1$, $a_2$ and $a_3$ be long vectors of
$V$ associated to the mirrors, extend them to an orthonormal
basis of $V$, and apply the lemma.
\proclaim
{Theorem}
{The variety $\calV$
is the intersection of the cubic eightfolds defined by the relations
of Lemma~\cubrel.}
Inter
\finishproclaim
{\it Proof sketch.}  Using one of the computer algebra systems {\pro
MACAULEY} or {\pro SINGULAR}, it is easy to see that the
dimension of the intersection $\calV'$ of the 270 cubics has
dimension 4.  With either system it is possible to compute a
projective resolution of $R/J$, where $R=\qz[Y_0,\dots,Y_9]$, 
$Y_0,\ldots,Y_9$ are a basis for $W$, and $J$ the
ideal generated by the 270 cubic relations.  The projective
dimension of $R/J$ turns out to be 5, by a 
calculation that takes a few minutes in {\pro SINGULAR} but several
hours in {\pro MACAULAY}.  As a consequence, $\calV'$
contains no component of dimension $<4$.
\smallskip
It is more involved to prove that $\calV'$ is irreducible. 
In principle one can simply ask the machine, but this seems to be too
much for the computer. Instead, we consider the intersection
of $\calV'$ with a hyperplane corresponding to a {cross}.
If $\calV'$ is irreducible then the intersection should consist
of~5 irreducible components. It is not hard to prove 
that in our situation 
the converse is also true. 
The hyperplane section is defined by a certain
ideal $\gota\subset\cz[Y_0,\dots,Y_9]$. In principle one can ask
the computer for the components of the ideal (e.g.\ by using
``decompose'' in {\pro MACAULEY}), but again this does not work.
Instead, one finds 
directly five ideals
$\gota_0,\dots,\gota_4$ containing $\gota$ that come from the
five mirrors of the {cross} and are
constructed in an obvious way. 
After the ideals $\gota_i$ have been constructed, one can verify
$\gota=\gota_0\cap\dots\cap\gota_4$ by means of {\pro MACAULEY}
or {\pro SINGULAR}.  The problem now is to prove that the
varieties of the $\gota_i$ are irreducible.  This means that
we face a similar problem in a lower dimension, which can be
treated in a similar manner. During this procedure several very
interesting ball quotients of smaller dimension occur. This will
be treated in a separate paper, where more details about the
ideal $J\subset\cz[Y_0,\dots,Y_9]$ and the hyperplane sections
will be given.
\smallskip
We also intend to include proofs of the facts that $J$ is
prime and that $R/J$ is normal. This has the important
consequence that $W$ generates the ring of all automorphic forms on
$\Gamma$ with trivial multipliers. The normality can be used to
give an alternate proof of Theorem~\embedding\ (that $\beta$ is
an embedding).  We will also give the
Hilbert function of $R/J$.
\qed
\neupara{Cross Ratios}%
In this section we will relate our automorphic forms to the
original invariants of a cubic surface, the cross-ratios of
Cayley. These are rational functions on the moduli space of
marked cubic surfaces that encode the manner in which the 27 lines on a
cubic surface lie in $P^3$. We will show below that Cayley's
cross-ratios are ratios of certain pairs of our 270
automorphic forms. Then we will use this to prove that the map
$\beta:\calM\to P^9$ of section~6
is an embedding.
\smallskip
Suppose that $A$ and $B$ are two crosses with the same long
mirror $m$. By Lemma~\EinF, the subspace of $W$ that changes sign
under reflection in $m$ is 2-dimensional, and the automorphic
forms coming from the six signed crosses of $m$ form a regular
hexagon in this plane, centered at $0$. Now, $A$ and $B$ define
two diameters of this hexagon, and we choose an endpoint $F$
(resp. $G$) of the diameter associated to $A$ (resp. $B$), such that
$F$ and $G$ are adjacent vertices of the hexagon. There is a
unique way to do this, up to simultaneously negating $F$ and
$G$, so the rational function $F/G$ does not depend on the
choice made. We call this the cross ratio $A/B$. The reason
for the name is Theorem~\ratios\  below, which identifies these
rational functions with Cayley's cross-ratios. It is a
genuine accident of terminology that our cross-ratios may
also be regarded as ratios of crosses. There are 270
cross-ratios, six for each of the 45 long mirrors. To identify
our cross-ratios with Cayley's we will need to describe the
divisor of $A/B$:
\proclaim
{Lemma}
{If $m$ is a long mirror in $X$ and $A$, $B$ and $C$ are its
three crosses, then the divisor of the cross-ratio $A/B$
consists of the four short mirrors of $A$ with multiplicity $1$
(simple zeros) and the four short mirrors of $B$ with
multiplicity $-1$ (simple poles). Furthermore, $A/B$  takes the
constant value $1$ at generic points of the short mirrors of
$C$.}
crdivisor%
\finishproclaim
\smallskip
{\it Proof.}  
If $F$ and $G$ are automorphic forms chosen as in
the discussion above, then their divisors in $X$ are the crosses
$A$ and $B$, respectively. Since the long mirrors of $A$ and $B$
coincide and the short mirrors are distinct, the identification
of the divisor of $A/B$ is complete. Finally,
$H=F-G$ is an endpoint of the third diameter of the hexagon, so that
it lies in the 1-dimensional subspace of $W$ associated to $C$,
and in particular it vanishes on the short mirrors of $C$. That
is, $F=G$ on the mirrors of $C$ and so $A/B=1$ at generic points
of the short mirrors of $C$.
\qed
\smallskip
Now we discuss Cayley's cross-ratios; our basic reference is
Naruki's
extensive study of them and a compactification $C$ of the
moduli space $M$ of marked smooth cubic surfaces that they
define [Na]. The biregular action of $W(E_6)$ on $M$ extends to a
biregular action on $C$, and the complement of $M$ in $C$ has 76
components, which fall into orbits of size 40 and 36 under
$W(E_6)$. The components in the orbit of size 40 are all
disjoint and can be blown down to points. The variety $\narukismodel$
obtained by this blowing-down is the standard Geometric
Invariant Theory (GIT) compactification of $M$, with its natural
$W(E_6)$-action. Now, $M$ is also $W(E_6)$-equivariantly
isomorphic to the complement in $X$ of the short mirrors, and
the inclusion of this space into $X$ is also the standard GIT
compactification. It follows that $X$ is
$W(E_6)$-equivariantly isomorphic to $\narukismodel$, with the 36 short
mirrors corresponding to the images in $\narukismodel$ of the remaining
36 components of $C-M$.
\smallskip
Naruki describes $M$ in terms of a maximal torus $T$ of the
simple Lie group $D_4$ of adjoint type. He writes $\Delta$ for
the union of the subtori which are the fixed-point sets of the
12 reflections of $W(D_4)$, and realizes $M$ as the blowup of
$T$ at the identity element, minus the proper transforms of the 12
components of $\Delta$. He introduces multiplicative characters $\lambda$,
$\mu$, $\nu$ and $\rho$ of $T$, which provide a coordinate
system for $T$, and describes the action of
$W(E_6)$ on $M$ by giving explicit rational self-maps of $T$ in terms
of these coordinates. This group $W(E_6)$ contains the obvious group
$W(D_4)$ and also the larger group $W(F_4)$ obtained by
adjoining the automorphisms of $T$ arising from the
automorphisms of the Dynkin diagram $D_4$. 
\smallskip
Naruki introduces 45
divisors in $M$ which $W(E_6)$ permutes transitively. One of
these, $\delta_0$, is the exceptional divisor lying
over the identity of $T$, and the rest are given by explicit
equations in $\lambda$, $\mu$, $\nu$ and $\rho$. We claim that
these 45 divisors correspond to our long mirrors. To see this,
observe that $\delta_0$ is the fixed-point set of the lift (say
$\eta$) to $M$ of the negation map on $T$. 
(All of the 2-torsion points of $T$ lie in $\Delta$.)
The conjugacy class
of $\eta$ has size at most 45, since $\eta$ centralizes
$W(F_4)$, of index 45 in $W(E_6)$, and at least 45, since
$\delta_0$ has 45 translates under $W(E_6)$. Since $W(E_6)$ has
a unique conjugacy class of involutions of size 45, and the
elements of this class are our long reflections, $\eta$
corresponds to a long reflection $\hat\eta$. Furthermore,
$\delta_0$ corresponds to the fixed-point set of $\hat\eta$,
which must therefore be irreducible and (since it contains a
long mirrors) consist entirely
of the long mirror. Another way to prove our claim is
to use the fact that each of [Na] and [ACT] proves that its set
of 45 divisors represent the marked cubic surfaces that have an
Eckardt point.
\smallskip
The passage from $T$ to $\narukismodel$ involves compactifying $T$ and
then performing a sequence of blowings-up and blowings-down. All
that matters to us is that the identity of $T$ is blown up, and
that the 12 components of $\Delta$ (or rather their
transforms in $\narukismodel$) are among the the 36
components of $\narukismodel-M$. Naruki calls these 12 divisors the
$A_1$-hypersurfaces.
\smallskip
Finally, Naruki's table~2 gives 45 of Cayley's cross-ratios
explicitly as rational functions of $\lambda$, $\mu$, $\nu$ and
$\rho$.  The full set of Cayley's 270 cross-ratios is obtained
by following these functions by the 6 projective linear
transformations of $P^1=\cz\cup\{\infty\}$ that preserve
$\{0,1,\infty\}$. Of course, Cayley had much more explicit
geometric concepts in mind when defining his cross-ratios; for
details see Naruki's paper.
\proclaim
{Theorem}
{Cayley's cross-ratios coincide with ours.}
ratios%
\finishproclaim
\smallskip
{\it Proof:}
The idea is to check that the divisors coincide and that
Cayley's cross-ratios satisfy the normalization condition
of lemma~\crdivisor. By Cayley's geometric considerations ([Na], \S3),
his cross-ratios do not take any of the values $0$, $1$ and $\infty$ in
$M$, so their divisors consist of short
mirrors with some multiplicities. For the short mirror $S$ given by
$\rho=1$ in Naruki's coordinates, simple substitution reveals
the behavior along $S$ of the 45 cross-ratios given in Naruki's
table. Namely, exactly 7 vanish along it, exactly 7 have poles
along it, and just one
takes the constant value 1.  Since the full set of Cayley's
cross-ratios is obtained by following these by the 6 linear
fractional transformations preserving $\{0,1,\infty\}$, we see
that exactly $2\cdot(7+7+1)=30$ of Cayley's cross-ratios vanish along
$S$, another 30 take the constant value 1, and a further 30
have poles along $S$. (Working with the full set of 270 restores
the symmetry between $0$, $1$ and $\infty$ that Naruki's choice
of 45 conceals.) Now, by the transitivity of $W(E_6)$ on
Cayley's cross-ratios, each vanishes along the same number, say
$k$, of short mirrors. By transitivity on the short mirrors,
each short mirror lies in the zero-locus of exactly 30 of
Cayley's cross-ratios. These transitivities also show that
$270\cdot k=36\cdot30$, so that $k=4$ and each of Cayley's
cross-ratios vanishes along exactly 4 short mirrors. The same
argument also shows that each has poles along exactly 4 short
mirrors. 
\smallskip
Now we consider Cayley's cross-ratio $r({\rm w})$, given in
Naruki's coordinates by
$$
r({\rm w})={
(\lambda\rho-1)(\mu\rho-1)(\nu\rho-1)(\lambda\mu\nu\rho-1)
\over
(\mu\nu\rho-1)(\lambda\nu\rho-1)(\lambda\mu\rho-1)(\rho-1)
}\;.
$$
We will write simply $r$ for $r({\rm w})$. The sets $\chi=1$,
where $\chi$ is one of the characters $\lambda\rho$, $\mu\rho$,
$\nu\rho$ and $\lambda\mu\nu\rho$ (resp. $\mu\nu\rho$,
$\lambda\nu\rho$, $\lambda\mu\rho$ and $\rho$) appearing in the numerator
(resp. denominator) are among Naruki's $A_1$-hypersurfaces, so
$r$ has a simple zero (resp. simple pole) along these four short
mirrors. By the argument above, these constitute the entire
divisor of $r$. Furthermore, the short mirrors along which
$r$ vanishes (resp. has a pole) are orthogonal, in the sense
that the reflections across them commute. To see this we do not even need to
perform a calculation, because Naruki (p.~20) has already
organized his twelve $A_1$-hypersurfaces into three sets each
consisting of four mutually orthogonal divisors. Finally, all 8
of these short mirrors are orthogonal to the long mirror
$\delta_0$, because their reflections obviously commute with the
negation map of $T$. It follows that $\delta_0$ together with
the four short mirrors coming from the numerator
(resp. denominator) of $r$ form a cross $C_n$
(resp. $C_d$). Therefore 
$r$ has the same divisor as our cross-ratio $C_n/C_d$. To
show that $r=C_n/C_d$ it now suffices to show that $r=1$ along
the short mirrors of the third cross associated to
$\delta_0$. Consulting again the table on Naruki's p.~20, we see
that these mirrors are given by $\chi=1$, where $\chi$ varies
over the characters 
$\lambda$, $\mu$, $\nu$ and $\lambda\mu\nu\rho^2$. Simple
calculation verifies the condition, so $r=C_n/C_d$. Since one
of Cayley's cross-ratios coincides with one of ours, and
$W(E_6)$ acts transitively on both sets of cross-ratios, the
theorem follows.
\qed
\smallskip
{\it Remark:}
B. van Geemen has also obtained this theorem, as a byproduct of
a larger investigation. His idea is to construct and study the
linear system on Naruki's model of the moduli space that comes
from our space $W$ of automorphic forms. After one understands
this linear system (van~Geemen identifies it with one introduced
by Coble long ago), the result above follows immediately.
His approach also has the advantage of allowing one to relate
the moduli space $\cal M$ to the variety $\cal V$ over fields
other than $\cz$. (Note that $\cal V$ is defined
over $\gz$.)
\proclaim
{Corollary}
{The map $\beta:X\to P^9$ of section~6 
is an embedding.}
embedding%
\finishproclaim
{\it Proof:}
We write $\partial X$ for
$\overline{\calB_4/\Gamma}-\calB_4/\Gamma$, the set of~40 cusp
points. One of Naruki's main results is that the~270
cross-ratios, a priori defined as maps
$M\to(P^1-\{0,1,\infty\})$, extend to regular maps $(X-\partial
X)\to P^1$ that embed $X-\partial X$ in $(P^1)^{270}$. Since the
cross-ratios are quotients of the elements of $W$, $\beta$ must
embed $X-\partial X$ in $P^9$. Unfortunately, this argument
cannot be extended to show that $\beta$ embeds all of $X$ into
$P^9$; the problem is that one must blow up the points of
$\partial X$ in order for the rational map from $X$ to $(P^1)^{270}$ to
become regular. In order to prove the theorem we will first show
that $\beta$ is injective as a map of sets, and then that
$\beta$ is a local embedding at each point of $\partial X$.
\smallskip
The injectivity has essentially already been proven: Theorem~\Fraemp\
shows that for each $x\in\partial X$, $x$ is the only point of
$X$ that lies on all the crosses containing $x$. It follows that
no point of $\partial X$ is identified under $\beta$ with any
other point of $X$. Since $\beta$ is already known to be
injective on $X-\partial X$, $\beta$ is injective.
\smallskip
Now we prove that $\beta$ is a local embedding at each point of
$\partial X$; we will use Naruki's explicit description (see [Na],
section~12) of these singularities. Namely, his
$T$-equivariant compactification $\tilde T$ of $T$ adjoins~48
divisors, 24 of which he then blows down to obtain~24 of the
points of $\partial X$. Focusing on one of these divisors, which
he denotes by $\overline{\rho=0}$ and we will denote by $D$, he
gives~8 characters of $T$ which extend to regular functions
$z_1,\ldots,z_8$ on a neighborhood $\calU$ of $D$ in $\tilde T$,
and which vanish along $D$. According to his theorem~12.1, the
induced map $\calU\to\cz^8$ gives the blowing-down of $D$ and thus
embeds a neighborhood of the resulting singular point
$x\in\partial X$ into $\cz^8$. Furthermore, he explicitly
describes the singularity as the cone on the Veronese embedding
of $P^1\times P^1\times P^1$ in $P^7$. This makes it a simple matter
to see that the divisor of each $z_i$ near $x$ has
exactly three components, and these components meet each
other away from $x$ as well as at $x$. Since each $z_i$ is the
extension of a character of $T$, its divisor can consists only
of the components of $X-M$, which is to say, short
mirrors. Since short mirrors that meet each other in $X-\partial
X$ must be orthogonal, we have shown that the divisor of each
$z_i$ near $X$ consists of three mutually orthogonal short
mirrors. For each $i=1,\ldots,8$, we will find an automorphic
form $\psi_i\in W$ whose divisor near $x$ coincides with that of
$z_i$. We may also choose $\psi'\in W$ whose divisor misses $x$
entirely. Then the evaluation of
$\psi_1/\psi',\ldots,\psi_8/\psi'$ provides essentially the same
map (of some neighborhood of $x$) into $\cz^8$ as Naruki's. It
follows that $\beta$ must embed a neighborhood of $x$ into
$P^9$.
\smallskip
All that remains is to show that if $x\in\partial X$ and $m_1$,
$m_2$ and $m_3$ are any three mutually orthogonal short mirrors
that all meet $x$, then there exists $\psi\in W$ whose divisor
near $X$ is just the sum of the $m_i$. We choose a primitive
null vector $v\in \Lambda$ representing $x$, and short roots
$r_i\in v^\perp$ whose mirrors represent the $m_i$. Denoting the
images of these vectors in $V$ by $\bar x$ and $\bar r_i$, we
may choose coordinates in $V$ so that the inner product is given
by 
$$
(a,b)=a_0b_0-a_1b_1-\cdots -a_4b_4\;,
$$
and $\bar r_1=(0,0,1,0,0)$, $\bar r_2=(0,0,0,1,0)$, $\bar
r_3=(0,0,0,0,1)$ and $\bar v=(1,1,0,0,0)$. The standard
cross, given by the pairs
$(\pm1,0,0,0,0),\ldots,(0,0,0,0,\pm1)$, is the divisor of one of
our Borcherds products, which we take to be $\psi$. It is
obvious that the divisor of $\psi$ contains the $m_i$. To show
that the other components of the divisor miss $x$, we observe
that these components correspond to the orthogonal complements
of the roots of $\Lambda$ whose images in $V$ are
$(\pm1,0,0,0,0)$ or $(0,\pm1,0,0,0)$. Any such root has inner
product $\not\equiv0$ (mod 3) with $v$, so its mirror cannot
contain $v$.
\qed
\vskip 1cm\noindent%
{\paragratit Literature}%
\bigskip
\def\leftheadline{\ninepoint\folio\hfill Literature}%
\def\rightheadline{\ninepoint Literature\hfill \folio}%
\item{[ACT1]} Allcock,\ D. Carlson,\ J. and Toledo, D.:
{\it A complex hyperbolic structure for the moduli of cubic
surfaces,}
Comptes Rendus de l'Aacademie Scientifique
Francaise {\bf 326}, ser~I, 49--54  (1988)\hfill\break
(alg geom /970916) 
\medskip
\item{[ACT2]} Allcock,\ D. Carlson,\ J., and Toledo, D.:
{\it The complex hyperbolic geometry of the moduli space
of cubic Surfaces,}  preprint
\medskip
\item{[BB]} Baily,\ W.L., and Borel,A.:
{\it Compactification of arithmetic quotients of bounded
symmetric domains,}  Annals of Math. {\bf 84},
No 3, 442--528  (1966)
\medskip
\item{[Bo1]} Borcherds,\ R.: {\it Automorphic forms
with singularities on Grassmannians,}
Invent. math. {\bf 132}, 491--562 (1998)
\hfill\break
http://www.dpmms.cam.ac.uk/home/emu/reb/.my-home-page.html
\medskip
\item{[Bo2]} Borcherds,\ R.: {\it The Gross-Kohnen-Zagier
theorem in higher dimensions, }
Duke Math.~J. {\bf 97}, No 219--233 (1999)
\hfill\break
http://www.dpmms.cam.ac.uk/home/emu/reb/.my-home-page.html
\medskip
\item{[Bo3]} Borcherds,\ R.: {\it An automorphic form related to
cubic surfaces,} 1997 (unpublished).
\medskip
\item{[C]} Conway, J. H. et. al.: {\it Atlas of Finite Groups,}
Oxford University Press, 1985.
\medskip
\item{[EHV]}  Eisenbud,\ D., Huneke,\ C., and Vasconcelos,\ W.:
{\it Direct methods for primary decomposition,}
Inv.~Math {\bf 110},  No 2,  207--235 (1992)
\medskip
\item{[Fr]} Freitag,\ E.:
{\it Modular forms related to cubic surfaces,}
in preperation.
\medskip
\item{[He]} Hecke,\ E.: Mathematische Werke,
G\"ottingen Vandenhoeck \&\ Ruprecht (1959)
\item{[Na]} Naruki, I.:
{\it Cross ratio variety as moduli space of cubic surfaces,}
Proc.~Lon\-don Math.~Soc. (3) {\bf 45}, 1--30  (1982)
\medskip
\item{[Sh]} Shimura,\ G.:
{\it The arithmetic of automorphic forms
with respect to a unitary group,}  
Annals of Math. {\bf 107}, 596--605 (1978)
\bye